\documentclass[12pt]{article}
\usepackage{amsmath,bbm,amssymb,mathrsfs}
\usepackage{amsthm,amscd,float}             
\usepackage{graphicx, multirow, graphics}
\usepackage{comment, color}
\usepackage{longtable, lscape}
\usepackage{rotating}
\usepackage{dcolumn}
\usepackage{natbib}
\usepackage{booktabs,caption,fixltx2e}
\usepackage[flushleft]{threeparttable}
\usepackage{caption}
\newtheorem{theorem}{Theorem}
\newtheorem{lemma}{Lemma}

\newtheorem{algorithm}{Algorithm}

\numberwithin{equation}{section}

\newcommand{\bm}[1]{\mbox{\boldmath$ #1 $\unboldmath}}

\def\tr{{\rm tr}}

\begin{document}
\baselineskip=22pt
\vskip 20pt

\begin{center}
{\Large \bf On Variable Ordination of Modified Cholesky Decomposition for Sparse Covariance Matrix Estimation}
\vskip 5pt
{Xiaoning Kang}$^{1}$ and {Xinwei Deng}$^{2}$\footnote{Address for correspondence: Xinwei Deng, Associate Professor, Department of Statistics,
Virginia Tech, Blacksburg, VA 24060 (E-mail: xdeng@vt.edu).}
\vskip 5pt
{\it{
$^{1}$International Business College,
Dongbei University of Finance and Economics, China \\
\vskip 3pt
$^{2}$Department of Statistics, Virginia Tech, Blacksburg, USA}} \\
\end{center}

\begin{abstract}
Estimation of large sparse covariance matrices is of great importance for statistical analysis, especially in the high-dimensional settings.
The traditional approach such as the sample covariance matrix performs poorly due to the high dimensionality.
The modified Cholesky decomposition (MCD) is a commonly used method for sparse covariance matrix estimation.
However, the MCD method relies on the order of variables, which is often not available or cannot be pre-determined in practice.
In this work, we solve this order issue by obtaining a set of covariance matrix estimates under different orders of variables used in the MCD.
Then we consider an ensemble estimator as the ``center'' of such a set of covariance matrix estimates with respect to the Frobenius norm.
The proposed method not only ensures the estimator to be positive definite, but also can capture the underlying sparse structure of the covariance matrix.
Under some weak regularity conditions, we establish both algorithmic convergence and asymptotical convergence of the proposed method.
The merits of the proposed method are illustrated through simulation studies and one real data example.
\end{abstract}


\section{Introduction}
Estimation of large covariance matrix from the high-dimensional data is an important and challenging problem in the multivariate data analysis.
For example, dimension reduction using the principal component analysis usually relies on accurate estimation of covariance matrix.
Under the context of graphical models, the estimation of covariance matrix or its inverse is often used to infer the network structure of the graph.
However, conventional estimation of covariance matrix is known to perform poorly due to the high dimensionality when the number of variables is close to or larger than the sample size (Johnstone, 2001).
To overcome the curse of dimensionality, a variety of methods proposed in literature often assumes certain patterns of sparsity for the covariance matrices.

In this work, our focus is on the estimation of sparse covariance matrix for high-dimensional data.
Early work on covariance matrix estimation includes shrinking eigenvalues of the sample covariance matrix (Dey \& Srinivasan, 1985; Haff, 1991),
a linear combination of the sample covariance and a proper diagonal matrix (Ledoit \& Wolf, 2004), improving the estimation based on matrix condition number (Aubry et al., 2012; Won et al., 2013),
and regularizing the eigenvectors of the matrix logarithm of the covariance matrix (Deng \& Tsui, 2013; Yu, Wang, \& Zhu, 2017).
However, the above mentioned methods do not explore the sparse structure of the covariance matrix.
A sparse covariance matrix estimate can be useful for subsequent inference, such as inferring the correlation pattern among the variables.
Bickel \& Levina (2009) proposed to threshold the small entries of the sample covariance matrix to zeroes and studied its theoretical behavior when the number of variables is large.
Rothman, Levina, \& Zhu (2009) considered to threshold the sample covariance matrix with more general thresholding functions.
Wagaman \& Levina (2009) introduced a method of the sparse estimation for the covariance matrix with banded structure based on the correlations between variables using the Isomap.
Cai \& Yuan (2012) proposed a covariance matrix estimation through block thresholding. Their estimator is constructed by dividing the sample covariance matrix into blocks and then simultaneously estimating the entries in a block by thresholding.
However, the threshold-based estimator is not guaranteed to be positive definite.
To make the estimate being sparse as well as positive definite, Bien \& Tibshirani (2011) considered a penalized likelihood method with a Lasso penalty (Tibshirani, 1996) on the entries of the covariance matrix.
Their idea is similar to the Graphical Lasso for inverse covariance matrix estimation in the literature (Yuan \& Lin, 2007; Friedman, Hastie, \& Tibshirani, 2008; Rocha, Zhao, \& Yu, 2008; Rothman et al., 2008; Yuan, 2008; Deng \& Yuan, 2009; and Yuan, 2010), but the computation is much more complicated due to the non-convexity of the objective function.
Xue, Ma, \& Zou (2012) developed a sparse covariance matrix estimator for high-dimensional data based on a convex objective function with positive definite constraint and $L_{1}$ penalty.
They also derived a fast algorithm to solve the constraint optimization problem.
Some other work on the estimation of high-dimensional covariance matrix can be found in Fan, Liao, \& Mincheva (2013), Liu, Wang, \& Zhao (2014), Xiao et al. (2016), Cai, Ren, \& Zhou (2016), Huang, Farewell, \& Pan (2017).
A comprehensive review of the development of covariance matrix estimation can be found in Pourahmadi (2013) and Fan, Liao, \& Liu (2016).

Another direction of sparse covariance matrix estimation is to take advantage of matrix decomposition.
One popular and effective decomposition is the modified Cholesky decomposition (MCD) (Pourahmadi, 1999; Wu \& Pourahmadi, 2003; Pourahmadi, Daniels, \& Park, 2007; Rothman, Levina, \& Zhu, 2009; Dellaportas \& Pourahmadi, 2012; Xue, Ma, \& Zou, 2012; Rajaratnam \& Salzman, 2013).
It assumes that the variables have a natural order based on which the variables can be sequentially orthogonalized to re-parameterize the covariance matrix.
By imposing certain sparse structures on the Cholesky factor, it results in certain sparse structure on the estimated covariance matrix.
For example,
Huang et al. (2006) considered to impose an $L_1$ (Lasso) penalty on the entries of the Cholesky factor for estimating the sparse covariance matrix.
Rothman, Levina, \& Zhu (2010) proposed a banded estimator of Cholesky factor, which can be obtained by regressing each variable only on its closest $k$ predecessors.
However, the MCD-based approach for estimating covariance matrix depends on the order of variables.
Such an pre-specification on the order of variables may not hold in practice.
A natural order of variables is often not available or cannot be pre-determined in many applications such as the gene expression data and stock marketing data.

In this paper, we adopt the MCD approach for estimating the large covariance matrix, but alleviate the drawback of order dependency of the MCD method using the permutation idea of Zheng et al. (2017).
By considering a set of covariance matrix estimates under different orders of variables in the MCD, Zheng et al. (2017) introduced an order-averaged estimator for the large covariance matrix with positive definite property.
However, such a approach cannot make the resultant estimator to be sparse.
In addition, they did not provide the theoretical results.
To overcome the drawbacks in Zheng et al. (2017), we address the order issue of MCD and encourage the sparse structure as well in the covariance matrix estimate.
We also show that the estimator of Zheng et al. (2017) can be considered as a special case of our proposed estimator when the penalty tuning parameter in the objective function is set to be zero.
It is worth to remarking that, under MCD for covariance matrix estimation, it is not straightforward to simultaneously address the order issue and sparsity together.
For this work, specifically, we first obtain a number of estimates of covariance matrix from different orders of variables using the permutation idea.
With such estimates, the proposed order-averaged estimator is obtained as the ``center'' of them under the Frobenius norm through an $L_{1}$ penalized objective function, where the $L_{1}$ regularization is imposed to achieve the sparsity of the estimate.
An efficient algorithm is also developed to make the computation attractive for solving the estimator.
Furthermore, the consistent property of the proposed estimator is also established under Frobenius norm with some regularity conditions.

The remainder of this work is organized as follows. Section \ref{sigma:sec:mcd} briefly reviews the MCD approach to estimate the covariance matrix. Section \ref{sigma:sec:methodology} introduces the proposed method by addressing the order issue.
An efficient algorithm is also developed to solve the objective function. In Section \ref{sigma:sec:converge}, the theoretical properties are presented. The simulation study and one real data example are reported in Section \ref{Sigma:simulation} and \ref{Sigma:application}, respectively. We conclude the paper in Section \ref{Sigma:sec:discussion}.

\section{Review of Modified Cholesky Decomposition}\label{sigma:sec:mcd}
Without loss of generality, suppose that $\bm X = (X_{1}, \ldots, X_{p})'$ is a $p$-dimensional random vector with mean $\bm 0$ and covariance matrix $\bm \Sigma$.
Let $\bm x_{1}, \ldots, \bm x_{n}$ be $n$ independent and identically distributed observations following $\mathcal{N}(\bm 0, \bm \Sigma)$.
Pourahmadi (1999) proposed the modified Cholesky decomposition (MCD) for the estimation of a covariance matrix, which is statistically meaningful and guarantees the positive definiteness of the estimate.
This decomposition arises from regressing each variable $X_{j}$ on its predecessors $X_{1}, \ldots, X_{j-1}$ for $2 \leq j \leq p$. Specifically, consider to fit a series of regressions
\begin{align*}
X_{j} = \sum_{k=1}^{j -1} (-t_{jk}) X_{k} + \epsilon_{j} = \hat{X}_{j} + \epsilon_{j},
\end{align*}
where $\epsilon_{j}$ is the error term for the $j$th regression with $E\epsilon_{j} = 0$ and $Var(\epsilon_{j}) = d_{j}^{2}$.
Let $\epsilon_{1} = X_{1}$ and $\bm D = diag(d_{1}^2, \ldots, d_{p}^2)$ be the diagonal covariance matrix of $\bm \epsilon = (\epsilon_{1}, \ldots, \epsilon_{p})'$. Construct the unit lower triangular matrix $\bm T = (t_{jk})_{p \times p}$ with ones on its diagonal and regression coefficients $(t_{j1}, \ldots, t_{j, j-1})'$ as its $j$th row. Then one can have
\begin{align*}
\bm D = Var(\bm \epsilon) = Var (\bm X - \hat{\bm X}) = Var (\bm T \bm X) = \bm T \bm \Sigma \bm T',
\end{align*}
and thus
\begin{align}\label{sigma:eq: mcd-1}
\bm \Sigma = \bm T^{-1} \bm D {\bm T'}^{-1}.
\end{align}
The MCD approach reduces the challenge of modeling a covariance matrix into the task of modeling $(p-1)$ regression problems, and is applicable in high dimensions.
However, directly imposing the sparse structure on Cholesky factor matrix $\bm T$ in \eqref{sigma:eq: mcd-1} does not imply the sparse pattern of covariance matrix $\bm \Sigma$
since it requires an inverse of $\bm T$. Thus the formulation \eqref{sigma:eq: mcd-1} is not convenient to impose a sparse structure on the estimation of $\bm \Sigma$.
Alternatively, one can consider a latent variable regression model based on the MCD. Writing $\bm X = \bm L \bm \epsilon$ would lead to
\begin{align}\label{sigma:MCD}
Var(\bm X) &= Var( \bm L \bm \epsilon )  \nonumber \\
\bm \Sigma &= \bm L \bm D \bm L'.
\end{align}
This decomposition can be interpreted as resulting from a new sequence of regressions, where each variable $X_{j}$ is regressed on all the previous latent variable $\epsilon_{1}, \ldots, \epsilon_{j-1} $ rather than themselves.
It gives a sequence of regressions
\begin{equation}\label{sigma:order}
X_j = \bm{l}^{T}_j \bm \epsilon = \sum_{k<j} l_{jk} \epsilon_{k} + \epsilon_j, \quad j=2, \ldots, p,
\end{equation}
where $\bm{l}_j=(l_{jk})$ is the $j$th row of $\bm L$. Here $l_{jj}=1$ and $l_{jk}=0$ for $k>j$.

With the data matrix $\mathbb{X} = (\bm x_{1}, \ldots, \bm x_{n})'$, define its $j$th column to be $\bm x^{(j)}$.
Denote by $\bm e^{(j)}$ the residuals of $\bm x^{(j)}$ for $j \geq 2$, and $\bm e^{(1)} = \bm x^{(1)}$.
Let $\mathbb{Z}^{(j)} = (\bm e^{(1)}, \ldots, \bm e^{(j-1)})$ be the matrix containing the first $(j-1)$ residuals.
Now the Lasso regularization (Tibshirani, 1996) can be used to encourage the sparsity on $\hat{\bm L}$ (Huang et al., 2006; Rothman, Levina, \& Zhu; 2010, Chang \& Tsay, 2010; Kang et al., 2019)
\begin{align}\label{sigma:eq:L}
\hat{\bm l}_{j} = \arg \min_{ \bm l_{j} } \| \bm x^{(j)} -  \mathbb{Z}^{(j)} \bm l_{j} \|_{2}^{2}
+ \eta_{j} \| \bm l_{j} \|_{1}, \ j = 2, \ldots, p,
\end{align}
where $\eta_{j} \geq 0$ is a tuning parameter and selected by cross validation. $\| \cdot \|_{1}$ stands for the vector $L_1$ norm. $\bm e^{(j)} = \bm x^{(j)} -  \mathbb{Z}^{(j)} \bm l_{j}$ is used to construct the residuals for the last column of $\mathbb{Z}^{(j+1)}$.
Then $d_{j}^2$ is estimated as the sample variance of $\bm e^{(j)}$
\begin{align}\label{sigma:eq:D}
\hat{d}_{j}^2 = \widehat{Var}(\hat{\bm e}^{(j)}) = \widehat{Var}(\bm x^{(j)} -  \mathbb{Z}^{(j)} \hat{\bm l}_{j})
\end{align}
when constructing matrix $\hat{\bm D} = diag(\hat{d}_{1}^2, \ldots, \hat{d}_{p}^2)$.
Hence, $\hat{\bm \Sigma} = \hat{\bm L} \hat{\bm D} \hat{\bm L'}$ will be a sparse covariance matrix estimate.

\section{The Proposed Method}\label{sigma:sec:methodology}

Clearly, the estimate $\hat{\bm \Sigma}= \hat{\bm L} \hat{\bm D} \hat{\bm L'}$ depends on the order of variables $X_{1}, \ldots, X_{p}$.
It means that different orders would lead to different estimates of $\bm \Sigma$.
To address this order-dependent issue, we consider an order-averaged estimation of $\bm \Sigma$ by using the idea of permutation.
Specifically, $M$ different permutations of $\{1, \ldots, p\}$ are generated as the orders of variables, denoted by
${\pi_{k}}^{,}s$, $k = 1, \ldots, M$. Let $\bm P_{\pi_{k}}$ be the corresponding permutation matrix.
Under a variable order $\pi_{k}$, the estimate is obtained as
\begin{align}\label{sigmapik}
\hat{\bm \Sigma}_{\pi_{k}} = \hat{\bm L}_{\pi_{k}} \hat{\bm D}_{\pi_{k}} \hat{\bm L}_{\pi_{k}}^{'},
\end{align}
where $\hat{\bm L}_{\pi_{k}}$ and $\hat{\bm D}_{\pi_{k}}$ are calculated based on \eqref{sigma:eq:L} and \eqref{sigma:eq:D}. Then transforming back to the original order, we have
\begin{align}\label{sigmak}
\hat{\bm\Sigma}_{k} = \bm P_{\pi_{k}} \hat{\bm \Sigma}_{\pi_{k}} \bm P_{\pi_{k}}'.
\end{align}
To obtain a proper estimator for $\bm \Sigma$, Zheng et al. (2017) proposed $\bar{\bm \Sigma} = \frac{1}{M} \sum_{k=1}^{M} \hat{\bm \Sigma}_{k}$.
However, such an estimate is clearly not sparse since the sparse structure in $\hat{\bm \Sigma}_{k}$ is destroyed by the average.

In order to simultaneously achieve the positive definiteness and sparsity for the estimator, we propose to consider
\begin{align}\label{sigma:equ1}
\hat{\bm \Sigma} = \arg \min_{ \bm \Sigma \succeq \nu \bm I} \frac{1}{2} \sum_{k=1}^{M} \| \bm \Sigma -  \hat{\bm \Sigma}_{k} \|_{F}^{2} + \tilde{\lambda} | \bm \Sigma |_{1},
\end{align}
where $\| \cdot \|_{F}$ stands for the Frobenius norm, $\tilde{\lambda} \geq 0$ is a tuning parameter, and $| \cdot |_{1}$ is $L_1$ norm for all the off-diagonal elements.
Here $\nu$ is some positive arbitrarily small number. The constraint $\bm \Sigma \succeq \nu \bm I$ is to guarantee the positive definiteness of the estimate.
The penalty term is to encourage the sparse pattern in $\hat{\bm \Sigma}$.
It is worth pointing out that, if $\tilde{\lambda} = 0$ in \eqref{sigma:equ1}, the solution of $\hat{\bm \Sigma}$ would be $\bar{\bm \Sigma} = \frac{1}{M} \sum_{k=1}^{M} \hat{\bm \Sigma}_{k}$, which is the estimator of Zheng et al. (2017);
if without the constraint $\bm \Sigma \succeq \nu \bm I$ in \eqref{sigma:equ1},
the solution of $\hat{\bm \Sigma}$ would be the soft-threshold estimate of $\bar{\bm \Sigma}$.
The objective \eqref{sigma:equ1} is similar to that of Xue, Ma, \& Zou (2012),
but their implications are different.
Xue, Ma, \& Zou (2012) used the sample covariance matrix $\bm S$ instead of
$\hat{\bm \Sigma}_{k}$ in \eqref{sigma:equ1}.
Hence, their estimate can be considered to have the minimum distance to $\bm S$
in terms of Frobenius norm.
However, our proposed estimate is to pursue the minimal averaged distance to all $\hat{\bm \Sigma}_{k}$'s, while maintain the properties of being positive definite and sparse.
The proposed estimate can be more accurate than the estimator proposed by Xue, Ma, \& Zou (2012), as evidenced in the numerical study of Section \ref{Sigma:simulation}.

For the ease of theoretical deduction, we re-write equation \eqref{sigma:equ1} as
\begin{align}\label{sigma:equ2}
\hat{\bm \Sigma} = \arg \min_{ \bm \Sigma \succeq \nu \bm I} \frac{1}{2 M} \sum_{k=1}^{M} \| \bm \Sigma -  \hat{\bm \Sigma}_{k} \|_{F}^{2} + \lambda | \bm \Sigma |_{1},
\end{align}
where $\lambda = \tilde{\lambda} / M$.
To efficiently solve the optimization \eqref{sigma:equ2}, we employ the alternating direction method of multipliers (ADMM) (Boyd et al., 2011), which has been widely used in solving the convex optimization of $L_1$ penalized covariance matrix estimation.
Let us first introduce a new variable $\bm \Phi$ and an equality constraint as follows
\begin{align}\label{sigma:equ3}
(\hat{\bm \Sigma}, \hat{\bm \Phi}) = \arg \min_{ \bm \Sigma, \bm \Phi }
\{\frac{1}{2 M} \sum_{k=1}^{M} \| \bm \Sigma -  \hat{\bm \Sigma}_{k} \|_{F}^{2}
+ \lambda | \bm \Sigma |_{1}: \bm \Sigma = \bm \Phi, \bm \Phi \succeq \nu \bm I \}.
\end{align}
Note that the solution of \eqref{sigma:equ3} gives solution to \eqref{sigma:equ2}.
To solve \eqref{sigma:equ3}, we minimize its augmented Lagrangian function for some given penalty parameter $\tau$ as
\begin{align}\label{sigma:optimization}
L(\bm \Sigma, \bm \Phi; \bm \Lambda) = \frac{1}{2 M} \sum_{k=1}^{M} \| \bm \Sigma -  \hat{\bm \Sigma}_{k} \|_{F}^{2}
+ \lambda | \bm \Sigma |_{1} - \langle \bm \Lambda, \bm \Phi - \bm \Sigma \rangle + \frac{1}{2\tau}\| \bm \Phi - \bm \Sigma \|_{F}^{2},
\end{align}
where $\bm \Lambda$ is the Lagrangian multiplier.
The notation $\langle \cdot, \cdot \rangle$ stands for the matrix inner product as $\langle \bm A, \bm B \rangle = \sum_{i,j} a_{ij} b_{ij}$, where $a_{ij}$ and $b_{ij}$ are the elements of matrices $\bm A$ and $\bm B$.
The ADMM iteratively solves the following steps sequentially for $i = 0, 1, 2,\ldots$ till convergence
\begin{align}\label{admm: theta}
\bm \Phi ~ \mbox{step}: \bm \Phi^{i+1} = \arg \min_{\bm \Phi \succeq \nu \bm I} L(\bm \Sigma^{i}, \bm \Phi; \bm \Lambda^{i})
\end{align}
\begin{align}\label{admm: sigma}
\bm \Sigma ~ \mbox{step}: \bm \Sigma^{i+1} = \arg \min_{\bm \Sigma} L(\bm \Sigma, \bm \Phi^{i+1}; \bm \Lambda^{i})
\end{align}
\begin{align}\label{admm: lambda}
\bm \Lambda ~ \mbox{step}: \bm \Lambda^{i+1} = \bm \Lambda^{i} - \frac{1}{\tau}(\bm \Phi^{i+1} - \bm \Sigma^{i+1}).
\end{align}
Assume the eigenvalue decomposition of a matrix $\bm Z$ is $\sum_{i=1}^{p}\lambda_{i}\bm \xi^{'} \bm \xi$, and
define $(\bm Z)_{+} = \sum_{i=1}^{p} \mbox{max}(\lambda_{i}, \nu) \bm \xi_{i}^{'} \bm \xi_{i}$.
Then we develop the closed form for \eqref{admm: theta} as
\begin{align*}
\frac{\partial L(\bm \Sigma^{i}, \bm \Phi; \bm \Lambda^{i})}{\partial \bm \Phi} &=
- \bm \Lambda^{i} + \frac{1}{\tau}(\bm \Phi - \bm \Sigma^{i}) \triangleq 0 \\
\bm \Phi &= \bm \Sigma^{i} + \tau \bm \Lambda^{i} \\
\bm \Phi^{i+1} &= (\bm \Sigma^{i} + \tau \bm \Lambda^{i})_{+}.
\end{align*}
Next, define an element-wise soft threshold for each entry $z_{ij}$ in matrix $\bm Z$ as
$\bm s(\bm Z, \delta) = \{ \bm s(z_{ij}, \delta) \}_{1 \leq i, j \leq p }$ with
\begin{align*}
\bm s(z_{ij}, \delta) = \mbox{sign}(z_{ij}) \mbox{max} ( |z_{ij}| - \delta, 0)I_{\{i \neq j \}} + z_{ij}I_{\{i = j \}}.
\end{align*}
Then the solution of \eqref{admm: sigma} is derived as
\begin{align*}
\frac{\partial L(\bm \Sigma, \bm \Phi^{i+1}; \bm \Lambda^{i})}{\partial \bm \Sigma} &= \frac{1}{M}\sum_{k=1}^{M}
(\bm \Sigma - \hat{\bm \Sigma}_{k}) + \bm \Lambda^{i} + \frac{1}{\tau}(\bm \Sigma - \bm \Phi^{i+1}) + \lambda \mbox{sign}^{\ast}(\bm \Sigma) \triangleq 0 \\
(\tau+1) \bm \Sigma &= \tau (\frac{1}{M}\sum_{k=1}^{M} \hat{\bm \Sigma}_{k} - \bm \Lambda^{i}) + \bm \Phi^{i+1} - \lambda \tau \mbox{sign}^{\ast}(\bm \Sigma) \\
\bm \Sigma^{i+1} &= \{ \bm s(\tau (\frac{1}{M}\sum_{k=1}^{M} \hat{\bm \Sigma}_{k} - \bm \Lambda^{i}) + \bm \Phi^{i+1}, \lambda \tau) \} / (\tau + 1),
\end{align*}
where $\mbox{sign}^{\ast}(\bm \Sigma)$ means $\mbox{sign}(\bm \Sigma)$ with the diagonal elements replaced by $\bm 0$ vector.
Algorithm \ref{sigma:alg_admm} summarizes the developed procedure for solving \eqref{sigma:equ2} by using the ADMM technique.
\begin{algorithm} \label{sigma:alg_admm}
~

\textbf{Step 1}: Input initial values $\bm \Sigma_{init}$, $\bm \Lambda_{init}$ and $\tau$.

\textbf{Step 2}: $\bm \Phi^{i+1} = (\bm \Sigma^{i} + \tau \bm \Lambda^{i})_{+}$.

\textbf{Step 3}: $\bm \Sigma^{i+1} = \{ \bm s(\tau (\frac{1}{M}\sum_{k=1}^{M} \hat{\bm \Sigma}_{k} - \bm \Lambda^{i}) + \bm \Phi^{i+1}, \lambda \tau) \} / (\tau + 1)$.

\textbf{Step 4}: $\bm \Lambda^{i+1} = \bm \Lambda^{i} - \frac{1}{\tau}(\bm \Phi^{i+1} - \bm \Sigma^{i+1})$.

\textbf{Step 5}: Repeat Step 2 - 4 till convergence.
\end{algorithm}
This algorithm converges fast and produces the optimal solution of $\arg \min L(\bm \Sigma, \bm \Phi; \bm \Lambda)$ in \eqref{sigma:optimization}.
In practice, the initial value $\bm \Sigma_{init}$ is set to be the estimate $\bar{\bm \Sigma} = \frac{1}{M} \sum_{k=1}^{M} \hat{\bm \Sigma}_{k}$.
The $\bm \Lambda_{init}$ is set to be zero matrix, and $\tau = 2$ as well as $\nu = 10^{-4}$.
The optimal value of tuning parameter $\lambda$ in \eqref{sigma:optimization} is chosen based on Bayesian information criterion (BIC) (Yuan \& Lin, 2007)
\begin{align*}
\mbox{BIC}(\lambda) = -\log|\hat{\bm \Sigma}^{-1}_{\lambda}| + \tr[\hat{\bm \Sigma}^{-1}_{\lambda} \bm S]
+ \frac{\log n}{n} \sum_{i \leq j} \hat{e}_{ij}(\lambda),
\end{align*}
where $\bm S$ is the sample covariance matrix, $\hat{\bm \Sigma}_{\lambda} = (\hat{\sigma}_{ij}^{(\lambda)})_{p \times p}$ indicates the  estimate of $\bm \Sigma$ obtained by applying our algorithm with tuning parameter $\lambda$.
$\hat{e}_{ij}(\lambda) = 0$ if $\hat{\sigma}_{ij}^{(\lambda)} = 0$, and $\hat{e}_{ij}(\lambda) = 1$ otherwise.


\section{Theoretical Convergence}\label{sigma:sec:converge}

In this section, Theorem 1 states that the sequence $(\bm \Sigma^{i}, \bm \Phi^{i}, \bm \Lambda^{i})$ generated
by Algorithm \ref{sigma:alg_admm} from any starting point numerically converges to an optimal minimizer
$(\hat{\bm \Sigma}^{+}, \hat{\bm \Phi}^{+}, \hat{\bm \Lambda}^{+})$ of \eqref{sigma:optimization}, where $\hat{\bm \Lambda}^{+}$ is the optimal dual variable.
Theorem 2 demonstrates the asymptotical consistent property of the proposed estimator under some weak regularity conditions.
The details of the proofs of Theorems 1 - 2 are in the Appendix.
To facilitate the presentation of the proofs, we first introduce some notations.
Define a $2p$ by $2p$ matrix $\bm J$ as
\begin{align*}
\bm J = \left (
\begin{array}{ccccc}
\tau \bm I_{p \times p} & 0  \\
0 & \tau^{-1} \bm I_{p \times p}
\end{array}
\right).
\end{align*}
Let the notation $\| \cdot \|_{J}^{2}$ be $\| \bm U \|_{J}^{2} = \langle \bm U, \bm J \bm U \rangle$
and $\langle \bm U, \bm V \rangle_{J} = \langle \bm U, \bm J \bm V \rangle$.
Let $\bm \Sigma_{0} = (\sigma^{0}_{ij})_{p \times p} = \bm L_{0} \bm D_{0} \bm L_{0}'$ be the true covariance matrix for the observations $\mathbb{X} = (x_{ij})_{n \times p}$, and define the number of nonzero off-diagonal elements of $\bm \Sigma_{0}$ as $s_{0}$. Denote the maximal true variance in $\bm \Sigma_{0}$ by $\sigma_{max}$.
Let $Z_{\pi_{k}} = \{(j, k): k < j, l_{0jk}^{(\pi_{k})} \neq 0 \}$ be the collection of nonzero elements in the lower triangular part of matrix $\bm L_{0\pi_{k}}$. Denote by $s_{1}$ the maximum of the cardinality of $Z_{\pi_{k}}$ for $k = 1, 2, \ldots, M$.
Now we present the following lemma and theory.

\begin{lemma}{Lemma 1.}{\label{1}}
Assume that $(\hat{\bm \Sigma}^{+}, \hat{\bm \Phi}^{+})$ is an optimal solution of \eqref{sigma:equ3} and
$\hat{\bm \Lambda}^{+}$ is the corresponding optimal dual variable with the equality constraint
$\bm \Sigma = \bm \Phi$, then the sequence $(\bm \Sigma^{i}, \bm \Phi^{i}, \bm \Lambda^{i})$
generated by Algorithm \ref{sigma:alg_admm} satisfies
\begin{align*}
\| \bm W^{+} - \bm W^{i} \|_{J}^{2} - \| \bm W^{+} - \bm W^{i+1} \|_{J}^{2} \geq \| \bm W^{i} - \bm W^{i+1} \|_{J}^{2},
\end{align*}
where $\bm W^{+} = (\hat{\bm \Lambda}^{+}, \hat{\bm \Sigma}^{+})'$ and $\bm W^{i} = (\bm \Lambda^{i}, \bm \Sigma^{i})'$.
\end{lemma}

\begin{theorem}{Theorem 1. (Algorithmic Convergence)}{\label{theory:convergence}}
Suppose $\bm x_{1}, \ldots, \bm x_{n}$ are $n$ independent and identically distributed observations from $\mathcal{N}_{p}(\bm 0, \bm \Sigma)$.
Then the sequence $(\bm \Sigma^{i}, \bm \Phi^{i}, \bm \Lambda^{i})$ generated by Algorithm \ref{sigma:alg_admm}
from any starting point converges to an optimal minimizer of the objective function in \eqref{sigma:optimization}.
\end{theorem}

Theorem 1 demonstrates the convergence of Algorithm \ref{sigma:alg_admm}.
It automatically indicates that the sequence $\bm \Sigma^{i}, i = 1, 2, \ldots$,
produced by Algorithm \ref{sigma:alg_admm} converges to an optimal solution of the objective \eqref{sigma:equ2}.
We prove Lemma 1 and Theorem 1 following the ideas of Xue, Ma \& Zou (2012) via the Karush-Kuhn-Tucker conditions (Karush, 1939; Kuhn \& Tucker, 1951).
The proofs are presented in the Appendix.

In order to achieve the asymptotical consistent property of the proposed estimator, one needs a basic assumption that there exists a constant $\theta > 1$ such that the singular values of the true covariance matrix are bounded as
\begin{align}\label{assumption1}
1/\theta < sv_{p}(\bm \Sigma_{0}) \leq sv_{1}(\bm \Sigma_{0}) < \theta,
\end{align}
where we use $sv_{1}(A), sv_{2}(A), \ldots, sv_{p}(A)$ to indicate the singular values of matrix $\bm A$ in a decreasing order. They are the square root of the eigenvalues of matrix $\bm A \bm A'$.
This assumption is also made in Rothman et al. (2008), Lam \& Fan (2009) and Guo et al. (2011). The assumption guarantees the positive definiteness property and makes inverting the covariance matrix meaningful. The following lemma and theorem give the asymptotical convergence
property of the proposed estimator under the Frobenius norm.

\begin{lemma}{Lemma 2.}{\label{sigma:lemma2_Jiang}}
Let $\bm \Sigma_{0} = \bm L_{0} \bm D_{0} \bm L_{0}^{'}$ be the MCD of the true covariance matrix.
If the singular values of $\bm \Sigma_{0}$ are bounded, there exist constants $\theta_{1}$ and $\theta_{2}$ such that $0 < \theta_{1} < sv_{p}(\bm \Sigma_{0}) \leq sv_{1}(\bm \Sigma_{0}) < \theta_{2} < \infty$, then there exist constants $h_{1}$ and $h_{2}$ such that
\begin{align*}
h_{1} < sv_{p}(\bm L_{0}) \leq sv_{1}(\bm L_{0}) < h_{2},
\end{align*}
and
\begin{align*}
h_{1} < sv_{p}(\bm D_{0}) \leq sv_{1}(\bm D_{0}) < h_{2}.
\end{align*}
\end{lemma}

\begin{lemma}{Lemma 3.}{\label{sigma:lemma2}}
Suppose $\bm x_{1}, \ldots, \bm x_{n}$ are $n$ independent and identically distributed observations from $\mathcal{N}_{p}(\bm 0, \bm \Sigma)$.
Let $\bm \Sigma_{0\pi_{k}} = \bm L_{0\pi_{k}} \bm D_{0\pi_{k}} \bm L_{0\pi_{k}}^{'}$ be the MCD of the true covariance matrix regarding a variable order $\pi_{k}$. Under \eqref{assumption1}, assume that the tuning parameters $\eta_{j}$ in \eqref{sigma:eq:L} satisfy $\sum_{j=1}^p \eta_{j} = O(\log(p) / n)$ and $(s_{1} + p) \log (p) = o(n)$,
then we have
\begin{align*}
\| \hat{\bm L}_{\pi_{k}} - \bm L_{0\pi_{k}} \|_{F} \stackrel{P}{\rightarrow} 0 ~~~ and ~~~
\| \hat{\bm D}_{\pi_{k}} - \bm D_{0\pi_{k}} \|_{F} \stackrel{P}{\rightarrow} 0.
\end{align*}
\end{lemma}


Lemma 3 demonstrates the asymptotical convergence of the Cholesky factor matrices
$\hat{\bm L}_{\pi_{k}}$ and $\hat{\bm D}_{\pi_{k}}$.
Based on this result, we can derive the theoretical property of
$\hat{\bm \Sigma}_{\pi_k}$ under variable order $\pi_{k}$, which is used to prove the following Theorem.

\begin{theorem}{Theorem 2. (Asymptotical Convergence)}{\label{theory:consistent1}}
Assume all the conditions in Lemma 3 hold, and $\lambda = o((s_{0} + p)^{-1/2})$.
Under the condition that for all $|t| \leq \rho$ and $1 \leq i \leq n, 1 \leq j \leq p$
\begin{align*}
E\{\exp (t x_{ij}^2)\} \leq K.
\end{align*}
For any $m > 0$, set $$\lambda = c_{0}^2 \frac{\log~p}{n} + c_{1} \left( \frac{\log~p}{n} \right)^{1/2},$$ where
$$c_{0} = \frac{1}{2} e K \rho^{1/2} + \rho^{-1/2} (m + 1)$$ and
$$c_{1} = 2 K (\rho^{-1} + \frac{1}{4} \rho \sigma_{max}^2)~\exp(\frac{1}{2} \rho \sigma_{max})~+~2\rho^{-1} (m + 2).$$
Then we have $$|| \hat{\bm \Sigma}^{+} - \bm \Sigma_{0} ||_{F} \stackrel{P}{\rightarrow} 0.$$
\end{theorem}

Theorem 2 demonstrates the asymptotically consistent properties of the proposed estimator with respect to the Frobenius norm under some regular conditions.
This together with Theorem 1 implies that the estimate obtained from Algorithm \ref{sigma:alg_admm} is consistent.
We would like to remark that the constraint $\bm \Sigma \succeq \nu \bm I$ could increase the computational cost from iterations in Algorithm \ref{sigma:alg_admm},
but it guarantees the proposed estimator to be positive definite.
Without this constraint, the solution of the optimization problem \eqref{sigma:equ1} would become the soft-threshold estimate of $\bar{\bm \Sigma}$.
Moreover, such a constraint helps to establish the convergence of Algorithm \ref{sigma:alg_admm} as well as the proposed estimator. Please refer to the proof for details.

\section{Simulation Study}\label{Sigma:simulation}
In this section, we conduct a comprehensive simulation study to evaluate the performance of the proposed method.
Suppose that data $\bm x_{1}, \ldots, \bm x_{n}$ are generated independently from the normal distribution $\mathcal{N}(\bm 0, \bm \Sigma)$.
Here we consider the following five covariance matrix structures.

$\textbf{Model}$ 1. $\bm \Sigma_{1}$ = MA(0.5, 0.3), where MA stands for ``moving average". The diagonal elements are 1 with the first sub-diagonal elements 0.5 and the seconde sub-diagonal elements 0.3.

$\textbf{Model}$ 2. $\bm \Sigma_{2}$ = AR(0.5), where AR stands for ``autoregressive". The conditional covariance between any two random variables
$X_{i}$ and $X_{j}$ is fixed to be $0.5^{|i - j|}$, $1 \leq i, j \leq p$.

$\textbf{Model}$ 3. $\bm \Sigma_{3}$ is generated by randomly permuting rows and corresponding columns of $\bm \Sigma_{1}$.

$\textbf{Model}$ 4. $\bm \Sigma_{4}$ is generated by randomly permuting rows and corresponding columns of $\bm \Sigma_{2}$.

$\textbf{Model}$ 5. $\bm \Sigma_{5} = \bm \Theta + \alpha \bm I$, where the diagonal elements of $\bm \Theta$ are zeroes and $\bm \Theta_{ij} = \bm \Theta_{ji} = b * Unif(-1, 1)$ for $i \neq j$, where $b$ is from the Bernoulli distribution with probability 0.15 equal 1. Each off-diagonal element of $\bm \Theta$ is generated independently.
The value of $\alpha$ is gradually increased to make sure that $\bm \Sigma_{5}$ is positive definite.

Note that $\textbf{Models}$ 1-2 consider the banded or nearly-banded structures for the covariance matrix.
While the covariance matrices of $\textbf{Models}$ 3-4 do not have structured sparsity due to the random permutations.
$\textbf{Model}$ 5 is similarly used in Bien \& Tibshirani (2011), which is a more general sparse matrix with no structure.
Hence from the perspective of sparse pattern, $\textbf{Model}$ 5 is most general case and $\textbf{Models}$ 1-2 are the least general cases.
For each case, we generate the data with three settings of different sample sizes and variable sizes: (1) $n = 50, p = 30$; (2) $n = 50, p = 50$ and (3) $n = 50, p = 100$.
For the implementation of the proposed method in this work, we choose $M = 100$ in the simulation where $p$ is at the order of 100.
We have tried $M = 10, 30, 50, 100$ and $150$ as the number of randomly selected permutations from all the possible $p!$ permutations.
The performances are seen to be marginally improved when $M$ is larger than 30.
Please refer to Kang \& Deng (2020) for a detailed discussion and justification on the choice of $M$. In practice, we would suggest to choose a relatively large $M$ to pursue the accuracy of the estimate if the computational resources are available.
Otherwise, a moderate value of $M$ is recommended to balance the accuracy and computation efficiency for the proposed method.

The performance of the proposed estimator is examined in comparison with several other approaches, which are divided into three classes.
The first class is the sample covariance matrix $\bm S$ that serves as the benchmark. The second class is composed of three methods that deal with the variable order used in the MCD, including the MCD-based method with BIC order selection (BIC) (Dellaportas \& Pourahmadi, 2012), the best permutation algorithm (BPA) (Rajaratnam \& Salzman, 2013) and the proposed method (Proposed).
The third class of competing methods consists of five approaches, including Bien and Tibshirani's estimate (BT) (Bien \& Tibshirani, 2011), Bickel and Levina's estimate (BL) (Bickel \& Levina, 2009), Xue, Ma and Zou's estimate (XMZ) (Xue, Ma, \& Zou, 2012), Wagaman and Levina's Isoband estimate (IB) (Wagaman \& Levina, 2009) and  Rothman et al.'s estimate (RLZ) (Rothman, Levina, \& Zhu, 2010).

To measure the accuracy of covariance matrix estimates $\hat{\bm \Sigma} = (\hat{\sigma}_{ij})_{p \times p}$ obtained from each approach, we consider the F norm, entropy loss (EN), $L_{1}$ norm and mean absolute error (MAE), defined as follows:
\begin{align*}
\mbox{F}     &= \sqrt{\sum_{i=1}^{p} \sum_{j=1}^{p} (\hat{\sigma}_{ij} - \sigma_{ij})^2}, \\
\mbox{EN}    &= \tr [\bm \Sigma^{-1} \hat{\bm \Sigma}] - \log |\bm \Sigma^{-1} \hat{\bm \Sigma}| - p, \\
L_{1}~\mbox{norm} &= \max_{ j } \sum_{i} | \hat{\sigma}_{ij} - \sigma_{ij} |, \\
\mbox{MAE} &= \frac{1}{p} \sum_{i=1}^{p} \sum_{j=1}^{p} |\hat{\sigma}_{ij} - \sigma_{ij}|.
\end{align*}
In addition, to gauge the performance of capturing sparse structure, we consider the false selection loss (FSL), which is the summation of false positive (FP) and false negative (FN).
Here we say a FP occurs if a nonzero element in the true matrix is incorrectly estimated as a zero. Similarly, a FN occurs if a zero element in the true matrix is incorrectly identified as a nonzero.
The FSL is computed in percentage as (FP + FN) / $p^2$.
For each loss function above, Tables \ref{sigma:table_model1} -  \ref{sigma:table_model2} and Tables \ref{sigma:table_model3} -  \ref{sigma:table_model5} in the Appendix report the averages of the performance measures and their corresponding standard errors in the parentheses over 100 replicates.
For each model, the two methods with lowest averages regarding each measure are shown in bold.
Dashed lines in the tables represent the corresponding values not available due to matrix singularity.

For a short summary of the numerical results, it shows that the proposed method generally provides a better estimation accuracy than other approaches in comparison.
It is able to accurately catch the underlying sparse structure of the covariance matrix.
Although the IB estimate gives good performance on estimation, it cannot guarantee the resultant estimator to be positive definite.
When the underlying covariance matrix is banded or tapered, the proposed method is not as good as the RLZ method.
The reason is that the RLZ method targets on the banded covariance matrix.
When the underlying structure of covariance matrix is more general without any specification, the proposed method still performs well.
As in the high-dimensional cases, the advantage of the proposed method is even more evident.

We first analyze the performance results and demonstrate the mechanism of several methods from the perspective of covariance structures using F norm as example.
Since the IB method assumes the true matrix has banded pattern after re-ordering the variables, it shows a good performance regarding F norm for $\textbf{Models}$ 1-4,
as $\textbf{Models}$ 1-2 are banded matrices and $\textbf{Models}$ 3-4 can have the banded structure after certain permutation of variables.
However, the IB method is inferior to the proposed method in $\textbf{Model}$ 5, since this model represents a general sparse covariance matrix with no possible banded structure, even if we permutate the variable order.
The RLZ method performs well for $\textbf{Models}$ 1-2 under F norm, as it is designated to estimate the banded or tapered matrices.
But it is not suitable for $\textbf{Models}$ 3-5.
In addition, we observe that the BPA approach produces relatively low F loss for $\textbf{Model}$ 4. The reason is that this method is able to recover the variable order for the AR model.
Therefore, although the BPA, IB and RLZ methods perform well for the banded or tapered matrices, they are inferior to the the proposed method
when the true covariance is a general matrix with no sparse pattern.

\begin{table}
\footnotesize
\begin{center}
\caption{The averages and standard errors of estimates for $\textbf{Model}$ 1.}
\label{sigma:table_model1}
\begin{tabular}{rrrrrrrrrrrrrrrr}
\hline\hline
&  &F &EN &$L_{1}$ &MAE &FSL (\%) \\\hline
\multirow {9}*{$p = 30$}
&$\bm S$&4.43 (0.05)      & 12.42 (0.08)    & 5.21 (0.07)     & 3.51 (0.03)     & 83.96 (0.01) \\
&BIC    &3.25 (0.03)      & 7.22 (0.08)     & 2.96 (0.04)     & 1.76 (0.02)     & 52.32 (0.55) \\
&BPA    &2.98 (0.03)      & \textbf{6.02 (0.09)} & 2.74 (0.05)    & 1.55 (0.02)     & 46.53 (0.68) \\
&BT     &4.74 (0.01)      & 7.78 (0.03)     & 2.14 (0.01)     & 1.85 (0.00)     & 6.92 (0.10) \\
&BL     &3.32 (0.05)      & -               & 2.33 (0.05)     & 1.17 (0.02) & 6.81 (0.14) \\
&XMZ    &3.35 (0.04)      & 10.61 (0.17)    & \textbf{1.88 (0.01)} & 1.25 (0.01) & 7.13 (0.18) \\
&IB     &\textbf{2.95 (0.05)}      & -               & 2.22 (0.05)    & \textbf{1.14 (0.03)}    & 8.86 (0.46) \\
&RLZ    &\textbf{2.90 (0.02)}      & 9.36 (0.07) & \textbf{1.55 (0.01)} & \textbf{1.04 (0.01)} & \textbf{6.22 (0.00)} \\
&Proposed   &3.26 (0.03)  & \textbf{7.10 (0.11)} & 1.92 (0.02) & 1.22 (0.01) & \textbf{6.75 (0.15)} \\
\midrule
\hline
\multirow {9}*{$p = 50$}
&$\bm S$&7.25 (0.05)      & -               & 8.26 (0.07)     & 5.75 (0.03)     & 90.19 (0.01) \\
&BIC    &4.51 (0.03)      & 15.77 (0.21)    & 3.85 (0.06)     & 1.98 (0.01)     & 43.30 (0.47) \\
&BPA    &4.30 (0.03)      & \textbf{12.98 (0.16)} & 3.62 (0.08)  & 1.84 (0.02)     & 41.49 (0.63) \\
&BT     &6.10 (0.07)      & 15.07 (0.28)    & 2.42 (0.02)     & 1.93 (0.02)     & 10.68 (0.46) \\
&BL     &4.67 (0.05)      & -               & 2.41 (0.05)     & 1.27 (0.01) & 4.80 (0.06) \\
&XMZ    &4.64 (0.05)      & 20.45 (0.29) & \textbf{1.98 (0.01)} & 1.36 (0.01) & 5.13 (0.07) \\
&IB     &\textbf{4.08 (0.05)}      & -               & 2.52 (0.04)    & \textbf{1.21 (0.02)}    & 6.05 (0.27) \\
&RLZ    &\textbf{3.79 (0.02)}      & 16.20 (0.07) & \textbf{1.63 (0.01)} & \textbf{1.06 (0.00)} & \textbf{3.84 (0.00)} \\
&Proposed   &4.58 (0.03)  & \textbf{13.68 (0.10)} & 2.02 (0.01) & 1.35 (0.01) & \textbf{4.36 (0.06)} \\
\midrule
\hline
\multirow {9}*{$p = 100$}
&$\bm S$&14.40 (0.06)     & -               & 16.15 (0.12)    & 11.43 (0.03)    & 95.01 (0.00) \\
&BIC    &6.87 (0.03)      & 42.56 (0.45)    & 5.26 (0.07)     & 2.24 (0.01)     & 32.75 (0.36) \\
&BPA    &6.74 (0.03)      & 35.68 (0.38)    & 5.29 (0.11)   & 2.20 (0.02)     & 33.60 (0.38) \\
&BT     &8.54 (0.13)      & \textbf{28.78 (0.33)} & 2.40 (0.03)  & 1.87 (0.02)     & 4.08 (0.33) \\
&BL     &7.19 (0.04)      & -               & 2.67 (0.06)     & 1.42 (0.01) & 2.83 (0.02) \\
&XMZ    &14.39 (0.06)     & 364.15 (0.18)   & 16.14 (0.12)    & 11.42 (0.03)    & 94.96 (0.01) \\
&IB     &\textbf{5.89 (0.05)}     & -               & 2.80 (0.05)    & \textbf{1.23 (0.01)}    & 2.71 (0.07) \\
&RLZ    &\textbf{5.41 (0.02)}      & 33.40 (0.12)    & \textbf{1.69 (0.01)} & \textbf{1.08 (0.00)} & \textbf{1.96 (0.00)} \\
&Proposed   &7.06 (0.02)   & \textbf{31.28 (0.14)} & \textbf{2.12 (0.01)} & 1.49 (0.00) & \textbf{2.38 (0.02)} \\
\midrule
\hline
\end{tabular}
\end{center}
\end{table}

\begin{table}
\footnotesize
\begin{center}
\caption{The averages and standard errors of estimates for $\textbf{Model}$ 2.}
\label{sigma:table_model2}
\begin{tabular}{rrrrrrrrrrrrrrr}
\hline\hline
&  &F &EN &$L_{1}$ &MAE &FSL (\%) \\\hline
\multirow {9}*{$p = 30$}
&$\bm S$&4.39 (0.04)      & 12.58 (0.08)    & 5.10 (0.07)     & 3.49 (0.03)     & 46.66 (0.01) \\
&BIC    &3.35 (0.03)      & 5.35 (0.09)     & 2.94 (0.04)     & 1.91 (0.01)     & 43.47 (0.30) \\
&BPA    &3.16 (0.03)      & 4.60 (0.07)     & 2.82 (0.04)     & 1.78 (0.01)     & 42.53 (0.31) \\
&BT     &4.70 (0.01)      & 5.40 (0.03)     & 2.52 (0.01)     & 2.19 (0.00)     & 43.82 (0.08) \\
&BL     &3.43 (0.04)      & -               & 2.61 (0.04) & 1.57 (0.01) & 43.70 (0.20) \\
&XMZ    &3.48 (0.04)      & 5.82 (0.11)     & \textbf{2.27 (0.01)} & 1.64 (0.01) & 42.24 (0.20) \\
&IB     &\textbf{3.02 (0.04)}      & -               & 2.61 (0.04)    & \textbf{1.56 (0.24)}    & \textbf{41.88 (0.45)} \\
&RLZ    &\textbf{2.76 (0.02)}      & \textbf{3.16 (0.03)} & \textbf{1.89 (0.01)} & \textbf{1.34 (0.01)} & 43.56 (0.00) \\
&Proposed   &3.47 (0.03)  & \textbf{4.09 (0.06)} & 2.30 (0.01) & 1.66 (0.01) & \textbf{41.82 (0.16)} \\
\midrule
\hline
\multirow {9}*{$p = 50$}
&$\bm S$&7.36 (0.05)      & -               & 8.54 (0.08)     & 5.84 (0.03)     & 65.59 (0.01) \\
&BIC    &4.57 (0.02)      & 11.08 (0.20)    & 3.93 (0.08)     & 2.18 (0.01)     & 42.28 (0.24) \\
&BPA    &4.38 (0.03)      & 9.17 (0.15)     & 3.70 (0.06)     & 2.08 (0.02)     & 41.73 (0.28) \\
&BT     &6.06 (0.05)      & 10.07 (0.18)    & 2.69 (0.02)     & 2.26 (0.01)     & 30.34 (0.13) \\
&BL     &4.73 (0.04)      & -               & 2.91 (0.05)     & 1.68 (0.01)     & 29.45 (0.07) \\
&XMZ    &4.73 (0.04)      & 11.26 (0.20)    & \textbf{2.36 (0.01)} & 1.76 (0.01) & 28.85 (0.07) \\
&IB     &\textbf{4.20 (0.05)}      & -               & 2.92 (0.04)    & \textbf{1.64 (0.02)}    & \textbf{27.99 (0.25)} \\
&RLZ    &\textbf{3.59 (0.01)}      & \textbf{5.48 (0.04)} & \textbf{1.96 (0.01)} & \textbf{1.38 (0.00)} &28.68 (0.00) \\
&Proposed   &4.70 (0.02)  & \textbf{7.22 (0.06)} & 2.40 (0.01) & 1.77 (0.01) & \textbf{28.60 (0.07)} \\
\midrule
\hline
\multirow {9}*{$p = 100$}
&$\bm S$&14.40 (0.07)     & -               & 16.04 (0.12)    & 11.43 (0.04)    & 81.86 (0.00) \\
&BIC    &6.92 (0.03)      & 29.70 (0.55)    & 5.32 (0.08)     & 2.47 (0.01)     & 33.77 (0.25) \\
&BPA    &6.78 (0.03)      & 23.65 (0.36)    & 5.16 (0.11)     & 2.44 (0.02)     & 34.41 (0.31) \\
&BT     &8.34 (0.13)      & 21.09 (0.36)    & 2.80 (0.03)     & 2.24 (0.02)     & 17.01 (0.23) \\
&BL     &7.18 (0.04)      & -               & 3.04 (0.05)     & 1.82 (0.01) & 16.07 (0.02) \\
&XMZ    &14.39 (0.07)     & 369.27 (0.19)   & 16.03 (0.12)    & 11.42 (0.04)    & 81.83 (0.00) \\
&IB     &\textbf{6.53 (0.05)}      & -               & 3.25 (0.05)    & \textbf{1.64 (0.01)}    & \textbf{15.22 (0.06)} \\
&RLZ    &\textbf{5.13 (0.02)}      & \textbf{11.03 (0.07)} & \textbf{2.05 (0.02)} & \textbf{1.41 (0.00)} & \textbf{15.12 (0.00)} \\
&Proposed   &7.11 (0.02)  & \textbf{16.29 (0.09)} & \textbf{2.49 (0.01)} & 1.90 (0.00) & 15.63 (0.02) \\
\midrule
\hline
\end{tabular}
\end{center}
\end{table}

Next, we provide some insights of the results through methods and other loss functions.
Tables \ref{sigma:table_model1} - \ref{sigma:table_model2} summarize the comparison results for $\textbf{Models}$ 1 and 2, respectively.
From the perspective of competing methods, the sample covariance matrix $\bm S$, serving as a benchmark approach,  does not give the sparse structure and performs poorly under all the loss measures.
The BIC and BPA in the second class of approaches provide sparse covariance matrix estimates compared with $\bm S$, but their false selection loss (FSL) are considerably larger than the proposed method.
Moreover, the proposed method greatly outperforms the BIC and BPA regarding $L_{1}$ and MAE for all settings of $p = 30, 50$ and 100.
Although the proposed method is comparable to the BIC and BPA methods under EN criterion when $p = 30$, it performs slightly better when $p = 50$ and much better in the case of $p = 100$.

For the BT method, the proposed method significantly outperforms it in capturing the sparse structure for the cases of $p=50$ and $p=100$.
Furthermore, the proposed method gives superior performance to the BT with respect to all the other loss criteria.
In comparison with the BL method, the proposed method performs similarly. It is known that the BL method is asymptotically optimal for sparse covariance matrix (Bickel \& Levina, 2009).
However, its estimated covariance matrix does not guarantee to be positive definite, which would result in matrix singularity in computing EN loss.
Compared with XMZ approach, the proposed method is superior or comparable with respect to all the loss measures in the settings of $p = 30$ and 50.
In the high-dimensional case when $p = 100$, the proposed method performs much better than the XMZ approach.
The IB method performs well regarding MAE and comparably to the proposed method under FSL. But it has the singularity issue.
Finally, it is seen that the proposed method is not as good as the RLZ approach for $\textbf{Models}$ 1-2 regarding these loss criteria.
This is not surprising, since the covariance matrices of $\textbf{Models}$ 1-2 are banded and tapered respectively, and the RLZ approach is designated to estimate such covariance matrix structures.

Tables \ref{sigma:table_model3} - \ref{sigma:table_model4} present the comparison results for $\textbf{Models}$ 3 and 4, respectively.
Different from $\textbf{Models}$ 1-2, the covariance matrices under $\textbf{Models}$ 3-4 are unstructured. This implies that the RLZ approach does not have the advantage.
Hence, it is clearly seen that the proposed method performs much better than the RLZ approach, especially at capturing the sparse structure and with respect to EN loss.
Generally, the proposed method provides superior performance to other approaches, with similar comparison results as described under $\textbf{Models}$ 1-2.
For the most general covariance matrix with no sparse pattern of $\textbf{Model}$ 5, the proposed method
performs even better as shown in Table \ref{sigma:table_model5} since it does not impose any assumptions on the sparse pattern of the underlying covariance matrix.
The proposed method gives substantial good performance when dimension is large.

\section{Application}\label{Sigma:application}
In this section, a real prostate cancer data set (Glaab et al., 2012) is used to evaluate the performance of the proposed method in comparison with other approaches described in Section \ref{Sigma:simulation}.
It contains two classes with 50 normal samples and 52 prostate cancer samples, and 2135 gene expression values recorded for each sample. Data are available online at http://ico2s.org/datasets/microarray.html.
Since it includes a large number of variables, the variable screening procedure is performed through two sample t-test.
Specifically, for each variable, t-test is conducted against the two classes of the prostate cancer data
such that the variables corresponding to large values of test statistics are ranked as significant variables.
Then the top 50 significant variables as group 1 and the top 50 nonsignificant variables as group 2 are selected for data analysis.
By doing so, there is supposed to be some correlations within each group of variables, but weak dependence between group 1 variables and group 2 variables (Rothman, Levina, \& Zhu, 2009; Xue, Ma, \& Zou, 2012).
Data are centered within each class and then used for the analysis.
In this section, to make each variable at the same scale, we focus on the correlation matrix rather than the covariance matrix.

\begin{figure}[h]
\begin{center}
\scalebox{0.5}[0.5]{\includegraphics{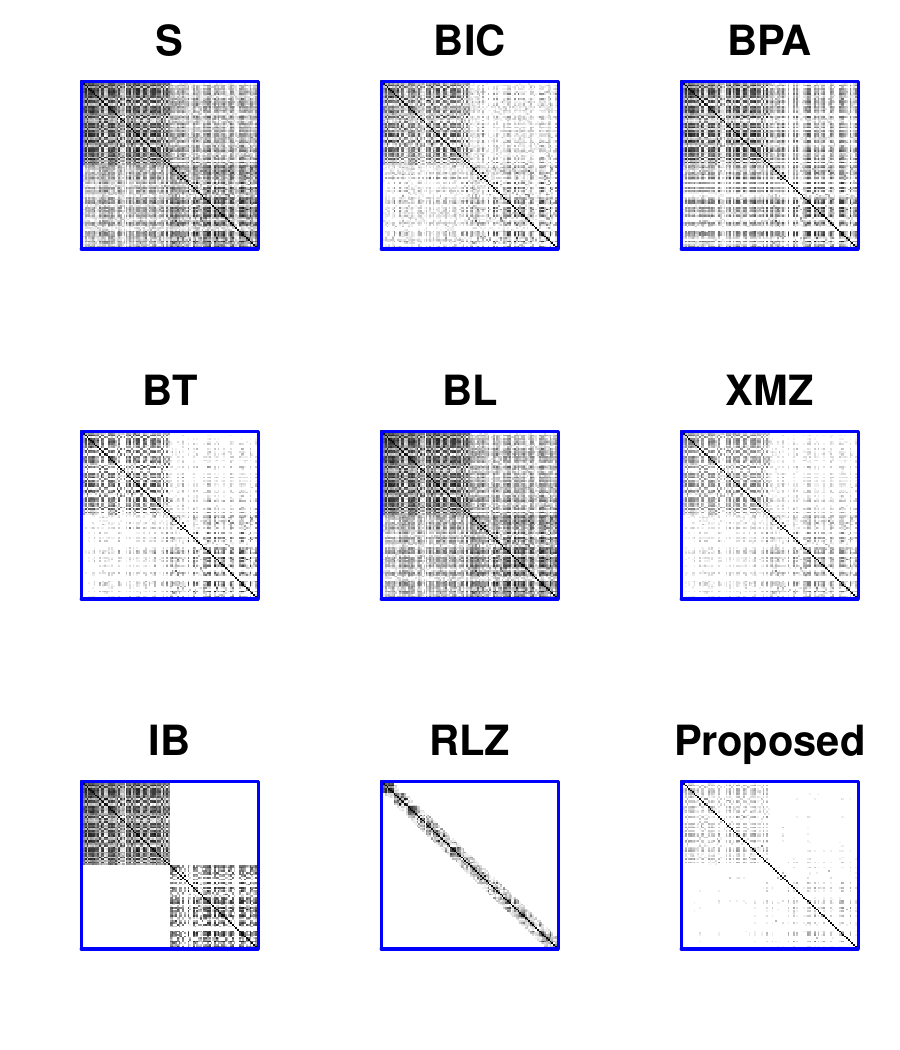}}
\caption{Heatmaps of the absolute values of the correlation matrices obtained from the proposed method and other approaches for prostate cancer data. Darker colour indicates higher density; lighter colour indicates lower density.} \label{figure:sigma_prostate}
\end{center}
\end{figure}

%
%

Figure \ref{figure:sigma_prostate} shows the heatmaps of the absolute values of the estimated correlation matrices obtained from each method.
It can be seen that for this data set, the IB method appears to have a leading performance for identifying the expected sparse pattern with clear blocks, followed by the proposed method, the BT and XMZ approaches, which are comparable to capture the sparse structure with two diagonal blocks.
All the rest approaches either result in a much sparser matrix as diagonal matrix (i.e. RLZ) or fail to identify the sparsity pattern (i.e. $\bm S$, BIC, BPA and BL).
We also observe that the IB and BL estimators yield negative eigenvalues, while the other estimators guarantee the positive definiteness.

\begin{table}[h]
\begin{center}
\caption{The misclassification errors in percentage of LDA.}\label{table:ME}
\begin{tabular}{lllllllllllll}
\hline\hline
&Methods    &BIC      &BPA      &BT       &XMZ     &RLZ     &Proposed \\
&ME         &31.7     &15.6     &14.9     &16.2    &14.4    &14.7     \\
&SD         &1.58     &1.06     &1.01     &0.97    &1.04    &0.95     \\
\hline\hline
\end{tabular}
\end{center}
\end{table}

Next, we further examine the performance of the proposed method by the classification of the linear discriminant analysis (LDA).
The whole data set is randomly split into the training set with 50 observations and the testing set with the rest 52 observations.
For this analysis, we screen all 2135 gene expressions by the two sample t-test based on the training data to select the top 100 significant variables.
Then all the compared methods use the training data to estimate the covariance matrix of these 100 variables.
Finally, each estimate is plugged into the LDA rule to classify the testing data.
Table \ref{table:ME} displays the misclassification errors (ME) in percentage and corresponding standard errors (SD) by each method for the above split procedure of 50 times. We see that although the proposed method is slightly inferior to the RLZ,
it performs better than other approaches in classification for this set of data.
Since we order the 100 variables by their significance from two sample t-test,
the variables far apart in distance from each other may have weak correlations.
Hence, the RLZ performs well as the covariance matrix of such 100 variables may be banded.

\section{Discussion}\label{Sigma:sec:discussion}
In this paper, we consider a positive definite estimate of covariance matrix based on the modified Cholesky decomposition (MCD).
The proposed method solves the order dependency issue in the MCD by considering the multiple estimates obtained from different variable orders.
The positive definite constraint and $L_{1}$ penalty are added to the objective function to guarantee the positive definiteness and encourage the sparse structure of the estimated covariance matrix. An efficient algorithm is developed to solve the constraint optimization problem. The proposed estimator does not require the prior knowledge of the variable order used in the MCD, and performs well in the high-dimensional cases. Simulation studies and real application demonstrate the superiority of the proposed method to several other existing approaches.

The idea of addressing variable ordination in this work may also be applied into other estimation problems, such as the inverse covariance matrix estimate.
However, one potential issue in practice is that the variables may have relations among themselves, i.e. causal relationship or spatial information.
It means that some orders of variables are meaningful and reflect such relations, while others not. This is can be clearly seen through BIC and BPA methods in the numerical study.
Hence, ruling out the meaningless orders and only using the meaningful orders of variables would improve the performance of the proposed method.
How to implement this idea in the real data needs further study.

\section*{BIBLIOGRAPHY}
\begin{description}
\item Aubry, A., De Maio, A., Pallotta, L., \& Farina, A. (2012). Maximum Likelihood Estimation of a Structured Covariance Matrix with a Condition Number Constraint.
\textit{Signal Processing, IEEE Transactions on}, 60(6), 3004--3021.

\item Bickel, P. J. \& Levina, E. (2004). Some Theory of Fisher's Linear Discriminant Function,Naive Bayes, and Some Alternatives When There Are Many More Variables than Observations.
\textit{Bernoulli}, 10(6), 989--1010.

\item Bickel, P. J. \& Levina, E. (2009). Covariance Regularization by Thresholding.
\textit{The Annals of Statistics}, 36(6), 2577--2604.

\item Bien, J. \& Tibshirani, R. J. (2011). Sparse Estimation of a Covariance Matrix.
\textit{Biometrika}, 98(4), 807--820.

\item Boyd, S., Parikh, N., Chu, E., Peleato, B., \& Eckstein, J. (2011). Distributed Optimization and Statistical Learning via the Alternating Direction Method of Multipliers.
\textit{Foundations and Trends in Machine Learning}, 3(1), 1--122.

\item Cai, T. T. \& Yuan, M. (2012). Adaptive Covariance Matrix Estimation through Block Thresholding.
\textit{The Annals of Statistics}, 40(40), 2014--2042.

\item Cai, T. T., Zhang, C. H., \& Zhou, H. H. (2010). Optimal Rates of Convergence for Covariance Matrix Estimation.
\textit{The Annals of Statistics}, 38(4), 2118--2144.

\item Cai, T. T., Ren, Z., \& Zhou, H. H. (2016). Estimating Structured High-Dimensional Covariance and Precision Matrices: Optimal Rates and Adaptive Estimation.
\textit{Electronic Journal of Statistics}, 10(1), 1--59.

\item Chang, C. \& Tsay, R. (2010). Estimation of Covariance Matrix via the Sparse Cholesky Factor with Lasso.
\textit{Journal of Statistical Planning and Inference}, 140(12), 3858--3873.

\item Dellaportas, P. \& Pourahmadi M. (2012). Cholesky-GARCH Models with Applications to Finance.
\textit{Statistics and Computing}, 22(4), 849--855.

\item Deng, X. \& Tsui, K. W. (2013). Penalized Covariance Matrix Estimation Using a Matrix-Logarithm Transformation.
\textit{Journal of Computational and Graphical Statistics}, 22(2), 494--512.

\item Deng, X. \& Yuan, M. (2009). Large Gaussian Covariance Matrix Estimation with Markov Structure,
\textit{Journal of Computational and Graphical Statistics},  18(3), 640--657.

\item Dey, D. K. \& Srinivasan, C. (1985). Estimation of a Covariance Matrix under Stein's Loss.
\textit{The Annals of Statistics}, 13(4), 1581--1591.

\item Fan, J., Liao, Y., \& Mincheva, M. (2013). Large Covariance Estimation by Thresholding Principal Orthogonal Complements.
\textit{Journal of the Royal Statistical Society: Series B}, 75(4), 603--680.

\item Fan, J., Liao, Y., \& Liu, H. (2016). An Overview of the Estimation of Large Covariance and Precision Matrices.
\textit{The Econometrics Journal}, 19(1), 1--32.

\item Friedman, J., Hastie, T., \& Tibshirani, T. (2008). Sparse Inverse Covariance Estimation with the Graphical Lasso.
\textit{Biostatistics}, 9(3), 432--441.

\item Glaab, E., Bacardit, J., Garibaldi, J. M., \& Krasnogor, N. (2012). Using Rule-Based Machine Learning for Candidate Disease Gene Prioritization and Sample Classification of Cancer Gene Expression Data.
\textit{PloS one}, 7(7), e39932.

\item Guo, J., Levina, E., Michailidis, G., \& Zhu, J. (2011). Joint Estimation of Multiple Graphical Models.
\textit{Biometrika},  98(1), 1--15.

\item Haff, L. R. (1991). The Variational Form of Certain Bayes Estimators.
\textit{The Annals of Statistics}, 19(3), 1163--1190.

\item Huang, C., Farewell, D., \& Pan, J. (2017). A Calibration Method for Non-positive Definite Covariance Matrix in Multivariate Data Analysis.
\textit{Journal of Multivariate Analysis}, 157, 45--52.

\item Huang, J. Z., Liu, N., Pourahmadi, M., \& Liu, L. (2006). Covariance Matrix Selection and Estimation via Penalised Normal Likelihood.
\textit{Biometrika},  93(1), 85--98.

\item Jiang, X. (2012). Joint Estimation of Covariance Matrix via Cholesky Decomposition.
Ph.D Dissertation. Department of Statistics and Applied Probability, National University of Singapore.

\item Johnstone, I. M. (2001). On the Distribution of the Largest Eigenvalue in Principal Components Analysis.
\textit{The Annals of Statistics}, 29(2), 295--327.

\item Kang, X., Deng, X., Tsui, K. W., \& Pourahmadi, M. (2019). On Variable Ordination of Modified Cholesky Decomposition for Estimating Time-Varying Covariance Matrices.
\textit{International Statistical Review}. DOI: 10.1111/insr.12357.

\item Kang, X. \& Deng, X. (2020). An Improved Modified Cholesky Decomposition Approach for Precision Matrix Estimation.
\textit{Journal of Statistical Computation and Simulation}, 90(3), 443--464.

\item Karush, W. (1939). Minima of Functions of Several Variables with Inequalities as Side Conditions.
Master Dissertation. Department of Mathematics, University of Chicago, Chicago, Illinois.

\item Kuhn, H. \& Tucker, A. (1951). Nonlinear Programming.
\textit{Proceedings of the 2nd Berkeley Symposium on Mathematica Statistics and Probabilistics},
University of California Press, 481--492.


\item Lam, C. \& Fan, J. (2009). Sparsistency and Rates of Convergence in Large Covariance Matrix Estimation.
\textit{The Annals of Statistics}, 37(6B), 4254--4278.

\item Lange, K., Hunter, D. R., \& Yang, I. (2000). Optimization Transfer Using Surrogate Objective Functions.
\textit{Journal of Computational and Graphical Statistics}, 9(1), 1--20.

\item Ledoit, O. \& Wolf, M. (2004). A Well-Conditioned Estimator for Large-Dimensional Covariance Matrices.
\textit{Journal of Multivariate Analysis}, 88(2), 365--411.


\item Liu, H., Wang, L., \& Zhao, T. (2014). Sparse Covariance Matrix Estimation with Eigenvalue Constraints.
\textit{Journal of Computational and Graphical Statistics}, 23(2), 439--459.

\item Pourahmadi, M. (1999). Joint Mean-Covariance Models with Applications to Longitudinal Data: Unconstrained Parameterisation.
\textit{Biometrika}, 86(3), 677--690.

\item Pourahmadi, M. (2013). High-Dimensional Covariance Estimation: with High-Dimensional Data. John Wiley \& Sons, Chichester, UK.

\item Pourahmadi, M., Daniels, M. J., \& Park, T. (2007). Simultaneous Modelling of the Cholesky Decomposition of Several Covariance Matrices.
\textit{Journal of Multivariate Analysis}, 98(3), 568--587.

\item Rajaratnam, B. \& Salzman, J. (2013). Best Permutation Analysis.
\textit{Journal of Multivariate Analysis},121(10), 193--223.

\item Rocha, G. V., Zhao, P., \& Yu, B. (2008). A Path Following Algorithm for Sparse Pseudo-Likelihood Inverse Covariance Estimation.
\textit{Technical Report}.

\item Rothman, A., Bickel, P., Levina, E., \& Zhu, J. (2008). Sparse Permutation Invariant Covariance Estimation.
\textit{Electronic Journal of Statistics}, 2(3), 494--515.

\item Rothman, A. J., Levina, E., \& Zhu, J. (2009). Generalized Thresholding of Large Covariance Matrices.
\textit{Journal of the American Statistical Association}, 104(485), 177--186.

\item Rothman, A. J., Levina, E., \& Zhu, J. (2010). A New Approach to Cholesky-Based Covariance Regularization in High Dimensions.
\textit{Biometrika}, 97(3), 539--550.

\item Tibshirani, R. (1996). Regression Shrinkage and Selection via the Lasso.
\textit{Journal of the Royal Statistical Society, Series B}, 58(1),  267--288.

\item Wagaman, A. \& Levina, E. (2009). Discovering Sparse Covariance Structures with the Isomap.
\textit{Journal of Computational and Graphical Statistics}, 18(3), 551--572.

\item Wu, W. B. \& Pourahmadi, M. (2003). Nonparametric Estimation of Large Covariance Matrices of Longitudinal Data.
\textit{Biometrika}, 90(4), 831--844.

\item Won, J. H., Lim, J., Kim, S. J., \& Rajaratnam, B. (2013). Condition - Number - Regularized Covariance Estimation.
\textit{Journal of the Royal Statistical Society: Series B}, 75(3), 427--450.

\item Xiao, L., Zipunnikov, V., Ruppert, D., \& Crainiceanu, C. (2016). Fast Covariance Estimation for High-Dimensional Functional Data.
\textit{Statistics and Computing}, 26(1-2), 409--421.

\item Xue, L., Ma, S., \& Zou, H. (2012). Positive-Definite $L_1$-Penalized Estimation of Large Covariance Matrices.
\textit{Journal of the American Statistical Association}, 107(500), 1480--1491.

\item Yu, P. L. H., Wang, X., \& Zhu, Y. (2017). High Dimensional Covariance Matrix Estimation by Penalizing the Matrix-Logarithm Transformed Likelihood.
\textit{Computational Statistics and Data Analysis}, 114, 12--25.

\item Yuan, M. \& Lin, Y. (2007). Model Selection and Estimation in the Gaussian Graphical Model.
\textit{Biometrika}, 94(1), 19--35.

\item Yuan, M. (2008). Efficient Computation of the $\ell_{1}$ Regularized Solution Path in Gaussian Graphical Models.
\textit{Journal of Computational and Graphical Statistics}, 17(4), 809--826.

\item Yuan, M. (2010). High Dimensional Inverse Covariance Matrix Estimation via Linear Programming.
\textit{The Journal of Machine Learning Research}, 11(12), 2261--2286.

\item Zheng H., Tsui K., Kang X., \& Deng X. (2017). Cholesky-Based Model Averaging for Covariance Matrix Estimation.
\textit{Statistical Theory and Related Fields}, 1(1), 48--58.
\end{description}

%

\begin{appendix}
\begin{proof}{Proof of Lemma 1}{}
Since $(\bm \Sigma^{+}, \bm \Phi^{+}, \bm \Lambda^{+})$ is the optimal
minimizer of \eqref{sigma:optimization}, based on the Karush-Kuhn-Tucker conditions we have
\begin{align}\label{equ:proof1}
(-\hat{\bm \Sigma}^{+} + \frac{1}{M}\sum_{k=1}^{M} \hat{\bm \Sigma}_{k} - \hat{\bm \Lambda}^{+})_{jl} \in \lambda \partial | \hat{\bm \Sigma}^{+}_{jl} |,
~j = 1,\ldots,p, l = 1,\ldots,p,  \mbox{and} ~ j \neq l
\end{align}
\begin{align}\label{equ:proof2}
(-\hat{\bm \Sigma}^{+} + \frac{1}{M}\sum_{k=1}^{M} \hat{\bm \Sigma}_{k})_{jj} + \hat{\bm \Lambda}^{+}_{jj} = 0, ~ j = 1,\ldots,p
\end{align}
\begin{align}\label{equ:proof3}
\hat{\bm \Phi}^{+} = \hat{\bm \Sigma}^{+}
\end{align}
\begin{align}\label{equ:proof4}
\hat{\bm \Phi}^{+} \succeq \nu \bm I,
\end{align}
and
\begin{align}\label{equ:proof5}
\langle \hat{\bm \Lambda}^{+}, \bm \Phi - \hat{\bm \Phi}^{+} \rangle \leq 0, ~ \forall \bm \Phi \succeq \nu \bm I.
\end{align}
The expressions in \eqref{equ:proof1} and \eqref{equ:proof2} result from the stationarity,
and \eqref{equ:proof3} and \eqref{equ:proof4} are valid because of the primal feasibility.
By the optimality conditions of the problem \eqref{admm: theta} with respect to $\bm \Phi$, we obtain
\begin{align*}
\langle \bm \Lambda^{i} - \frac{1}{\tau}(\bm \Phi^{i+1} - \bm \Sigma^{i}), \bm \Phi - \bm \Phi^{i+1} \rangle \leq 0,
~ \forall \bm \Phi \succeq \nu \bm I.
\end{align*}
This, together with $\bm \Lambda$ step \eqref{admm: lambda}, yields
\begin{align}\label{equ:proof6}
\langle \bm \Lambda^{i+1} - \frac{1}{\tau}(\bm \Sigma^{i+1} - \bm \Sigma^{i}), \bm \Phi - \bm \Phi^{i+1} \rangle \leq 0,
~ \forall \bm \Phi \succeq \nu \bm I.
\end{align}
Now by setting $\bm \Phi = \bm \Phi^{i+1}$ in \eqref{equ:proof5} and $\bm \Phi = \hat{\bm \Phi}^{+}$ in \eqref{equ:proof6} respectively, it leads to
\begin{align}\label{equ:proof7}
\langle \hat{\bm \Lambda}^{+}, \bm \Phi^{i+1} - \hat{\bm \Phi}^{+} \rangle \leq 0,
\end{align}
and
\begin{align}\label{equ:proof8}
\langle \bm \Lambda^{i+1} - \frac{1}{\tau}(\bm \Sigma^{i+1} - \bm \Sigma^{i}), \hat{\bm \Phi}^{+} - \bm \Phi^{i+1} \rangle \leq 0.
\end{align}
Summing \eqref{equ:proof7} and \eqref{equ:proof8} gives
\begin{align}\label{equ:proof9}
\langle (\bm \Lambda^{i+1} - \hat{\bm \Lambda}^{+}) - \frac{1}{\tau}(\bm \Sigma^{i+1} - \bm \Sigma^{i}), \bm \Phi^{i+1} - \hat{\bm \Phi}^{+} \rangle \geq 0.
\end{align}
On the other hand, by the optimality conditions of the problem \eqref{admm: sigma} with respect to $\bm \Sigma$, we have
\begin{align}\label{equ:proof10}
0 \in [ \frac{1}{M}\sum_{k=1}^{M}(\bm \Sigma^{i+1} - \hat{\bm \Sigma}_{k}) + \bm \Lambda^{i} + \frac{1}{\tau}(\bm \Sigma^{i+1} - \bm \Phi^{i+1})]_{jl}
+ \lambda \partial |\bm \Sigma^{i+1}_{jl}|, ~  j \neq l,
\end{align}
and
\begin{align}\label{equ:proof11}
[ \frac{1}{M}\sum_{k=1}^{M}(\bm \Sigma^{i+1} - \hat{\bm \Sigma}_{k}) + \bm \Lambda^{i} + \frac{1}{\tau}(\bm \Sigma^{i+1} - \bm \Phi^{i+1})]_{jj}
= 0, ~ j = 1,\ldots,p,
\end{align}
Plugging $\bm \Lambda$ step \eqref{admm: lambda} into \eqref{equ:proof10} and \eqref{equ:proof11} respectively
results in
\begin{align}\label{equ:proof12}
( - \bm \Sigma^{i+1} + \frac{1}{M}\sum_{k=1}^{M}\hat{\bm \Sigma}_{k} - \bm \Lambda^{i+1})_{jl} \in
\lambda \partial |\bm \Sigma^{i+1}_{jl}|, ~ j = 1,\ldots,p, l = 1,\ldots,p,  \mbox{and} ~ j \neq l,
\end{align}
and
\begin{align}\label{equ:proof13}
(\bm \Sigma^{i+1} - \frac{1}{M}\sum_{k=1}^{M}\hat{\bm \Sigma}_{k})_{jj} + \bm \Lambda^{i+1}_{jj} = 0, ~ j = 1,\ldots,p.
\end{align}
Since $\partial | \cdot |$ is monotonically non-decreasing, \eqref{equ:proof1} and \eqref{equ:proof12} yield for $j \neq l$
\begin{align*}
( - \bm \Sigma^{i+1} + \frac{1}{M}\sum_{k=1}^{M}\hat{\bm \Sigma}_{k} - \bm \Lambda^{i+1})_{jl} \left\{
\begin{array}{l}
\geq ( - \hat{\bm \Sigma}^{+} + \frac{1}{M}\sum_{k=1}^{M} \hat{\bm \Sigma}_{k} - \hat{\bm \Lambda}^{+})_{jl}, ~\mbox{if}~~\bm \Sigma^{i+1}_{jl} \geq \hat{\bm \Sigma}^{+}_{jl} \\
\leq ( - \hat{\bm \Sigma}^{+} + \frac{1}{M}\sum_{k=1}^{M} \hat{\bm \Sigma}_{k} - \hat{\bm \Lambda}^{+})_{jl}, ~\mbox{if}~~\bm \Sigma^{i+1}_{jl} < \hat{\bm \Sigma}^{+}_{jl}
\end{array},
\right.
\end{align*}
that is,
\begin{align*}
(\hat{\bm \Sigma}^{+} - \bm \Sigma^{i+1} + \hat{\bm \Lambda}^{+} - \bm \Lambda^{i+1})_{jl} \left\{
\begin{array}{l}
\geq 0, ~\mbox{if}~~\bm \Sigma^{i+1}_{jl} \geq \hat{\bm \Sigma}^{+}_{jl} \\
\leq 0, ~\mbox{if}~~\bm \Sigma^{i+1}_{jl} < \hat{\bm \Sigma}^{+}_{jl}
\end{array}.
\right.
\end{align*}
As a result, we obtain
\begin{align}\label{equ:proof14}
(\bm \Sigma^{i+1} - \hat{\bm \Sigma}^{+})_{jl}(\hat{\bm \Sigma}^{+} - \bm \Sigma^{i+1} + \hat{\bm \Lambda}^{+} - \bm \Lambda^{i+1})_{jl} \geq 0,
~ j = 1,\ldots,p, l = 1,\ldots,p,  \mbox{and} ~ j \neq l.
\end{align}
In addition, subtracting \eqref{equ:proof13} from \eqref{equ:proof2} implies
\begin{align}\label{equ:proof15}
(\hat{\bm \Sigma}^{+} - \bm \Sigma^{i+1} + \hat{\bm \Lambda}^{+} - \bm \Lambda^{i+1})_{jj} = 0, ~ j = 1, \ldots, p.
\end{align}
Then combining \eqref{equ:proof14} and \eqref{equ:proof15} leads to
\begin{align}\label{equ:proof16}
\langle \bm \Sigma^{i+1} - \hat{\bm \Sigma}^{+}, \hat{\bm \Sigma}^{+} - \bm \Sigma^{i+1} + \hat{\bm \Lambda}^{+} - \bm \Lambda^{i+1} \rangle \geq 0.
\end{align}
By summing \eqref{equ:proof9} and \eqref{equ:proof16}, we have
\begin{align*}
\langle \bm \Sigma^{i+1} - \hat{\bm \Sigma}^{+}, \hat{\bm \Lambda}^{+} - \bm \Lambda^{i+1} \rangle +
\langle \bm \Lambda^{i+1} - \hat{\bm \Lambda}^{+}, \bm \Phi^{i+1} - \hat{\bm \Phi}^{+} \rangle -
\frac{1}{\tau} \langle \bm \Sigma^{i+1} - \hat{\bm \Sigma}^{i}, \bm \Phi^{i+1} - \hat{\bm \Phi}^{+} \rangle \\
\geq \| \bm \Sigma^{i+1} - \hat{\bm \Sigma}^{+} \|^{2}_{F}.
\end{align*}
This, together with \eqref{equ:proof3} and $\bm \Phi^{i+1} = \tau (\bm \Lambda^{i} - \bm \Lambda^{i+1}) + \bm \Sigma^{i+1}$
from $\bm \Lambda$ step \eqref{admm: lambda}, gives
\begin{align}\label{equ:proof17}
\tau \langle \bm \Lambda^{i+1} - \hat{\bm \Lambda}^{+}, \bm \Lambda^{i} - \bm \Lambda^{i+1} \rangle +
\frac{1}{\tau} \langle \bm \Sigma^{i+1} - \hat{\bm \Sigma}^{+}, \bm \Sigma^{i} - \bm \Sigma^{i+1} \rangle \nonumber \\
\geq \| \bm \Sigma^{i+1} - \hat{\bm \Sigma}^{+} \|^{2}_{F} - \langle \bm \Lambda^{i} - \bm \Lambda^{i+1}, \bm \Sigma^{i} - \bm \Sigma^{i+1} \rangle
\end{align}
By $\hat{\bm \Phi}^{+} - \bm \Phi^{i+1} = (\hat{\bm \Phi}^{+} - \bm \Phi^{i}) + (\bm \Phi^{i} - \bm \Phi^{i+1})$ and
$\hat{\bm \Sigma}^{+} - \bm \Sigma^{i+1} = (\hat{\bm \Sigma}^{+} - \bm \Sigma^{i}) + (\bm \Sigma^{i} - \bm \Sigma^{i+1})$,
\eqref{equ:proof17} is reduced to
\begin{align}\label{equ:proof18}
\tau \langle \bm \Lambda^{i} - \hat{\bm \Lambda}^{+}, \bm \Lambda^{i} - \bm \Lambda^{i+1} \rangle +
\frac{1}{\tau} \langle \bm \Sigma^{i} - \hat{\bm \Sigma}^{+}, \bm \Sigma^{i} - \bm \Sigma^{i+1} \rangle
\geq \tau \| \bm \Lambda^{i} - \bm \Lambda^{i+1} \|^{2}_{F} \nonumber \\
+ \frac{1}{\tau} \| \bm \Sigma^{i} - \bm \Sigma^{i+1} \|^{2}_{F}
+ \| \bm \Sigma^{i+1} - \hat{\bm \Sigma}^{+} \|^{2}_{F} - \langle \bm \Lambda^{i} - \bm \Lambda^{i+1}, \bm \Sigma^{i} - \bm \Sigma^{i+1} \rangle
\end{align}
Using the notations $\bm W^{+}$ and $\bm W^{i}$, the left hand side of \eqref{equ:proof18} becomes
\begin{align*}
&~~~~ \langle (\bm \Lambda^{i} - \hat{\bm \Lambda}^{+}, \bm \Sigma^{i} - \hat{\bm \Sigma}^{+})',
[\tau (\bm \Lambda^{i} - \bm \Lambda^{i+1}), \frac{1}{\tau} (\bm \Sigma^{i} - \bm \Sigma^{i+1})]' \rangle \\
&= \langle (\bm \Lambda^{i}, \bm \Sigma^{i})' - (\hat{\bm \Lambda}^{+}, \hat{\bm \Sigma}^{+})',
\bm J [(\bm \Lambda^{i}, \bm \Sigma^{i})' - (\hat{\bm \Lambda}^{i+1}, \hat{\bm \Sigma}^{i+1})'] \rangle \\
&= \langle \bm W^{i} - \bm W^{+}, \bm J (\bm W^{i} - \bm W^{i+1}) \rangle \\
&= \langle \bm W^{i} - \bm W^{+}, \bm W^{i} - \bm W^{i+1} \rangle_{J}.
\end{align*}
The first two terms on the right side of \eqref{equ:proof18} becomes
\begin{align*}
\tau \| \bm \Lambda^{i} - \bm \Lambda^{i+1} \|^{2}_{F} + \frac{1}{\tau} \| \bm \Sigma^{i} - \bm \Sigma^{i+1} \|^{2}_{F}
&= \tau \langle \bm \Lambda^{i} - \bm \Lambda^{i+1}, \bm \Lambda^{i} - \bm \Lambda^{i+1} \rangle +
\frac{1}{\tau} \langle \bm \Sigma^{i} - \bm \Sigma^{i+1}, \bm \Sigma^{i} - \bm \Sigma^{i+1} \rangle \\
&= \langle (\bm \Lambda^{i} - \bm \Lambda^{i+1},  \bm \Sigma^{i} - \bm \Sigma^{i+1})',
[\tau (\bm \Lambda^{i} - \bm \Lambda^{i+1}), \frac{1}{\tau} (\bm \Sigma^{i} - \bm \Sigma^{i+1})]' \rangle \\
&= \langle (\bm \Lambda^{i}, \bm \Sigma^{i})' - (\bm \Lambda^{i+1}, \bm \Sigma^{i+1})',
\bm J [(\bm \Lambda^{i}, \bm \Sigma^{i})' - (\bm \Lambda^{i+1}, \bm \Sigma^{i+1})'] \rangle \\
&= \langle \bm W^{i} - \bm W^{i+1}, \bm J (\bm W^{i} - \bm W^{i+1}) \rangle \\
&= \| \bm W^{i} - \bm W^{i+1} \|^{2}_{J}.
\end{align*}
As a result, \eqref{equ:proof18} can be rewritten as
\begin{align}\label{equ:proof19}
\langle \bm W^{i} - \bm W^{+}, \bm W^{i} - \bm W^{i+1} \rangle_{J} \geq \| \bm W^{i} - \bm W^{i+1} \|^{2}_{J}
+ \| \bm \Sigma^{i+1} - \hat{\bm \Sigma}^{+} \|^{2}_{F} - \langle \bm \Lambda^{i} - \bm \Lambda^{i+1}, \bm \Sigma^{i} - \bm \Sigma^{i+1} \rangle.
\end{align}
Note a fact that
\begin{align*}
\| \bm W^{+} - \bm W^{i+1} \|^{2}_{J} = \| \bm W^{+} - \bm W^{i} \|^{2}_{J} - 2 \langle \bm W^{+} - \bm W^{i}, \bm W^{i+1} - \bm W^{i} \rangle_{J}
+ \| \bm W^{i} - \bm W^{i+1} \|^{2}_{J}.
\end{align*}
Therefore,
\begin{align}\label{equ:proof20}
&~~~~ \| \bm W^{+} - \bm W^{i} \|^{2}_{J} - \| \bm W^{+} - \bm W^{i+1} \|^{2}_{J}  \nonumber \\
&= 2 \langle \bm W^{+} - \bm W^{i}, \bm W^{i+1} - \bm W^{i} \rangle_{J} - \| \bm W^{i} - \bm W^{i+1} \|^{2}_{J} \nonumber \\
&\geq 2 \| \bm W^{i} - \bm W^{i+1} \|^{2}_{J} + 2 \| \bm \Sigma^{i+1} - \hat{\bm \Sigma}^{+} \|^{2}_{F}
- 2 \langle \bm \Lambda^{i} - \bm \Lambda^{i+1}, \bm \Sigma^{i} - \bm \Sigma^{i+1} \rangle - \| \bm W^{i} - \bm W^{i+1} \|^{2}_{J} \nonumber \\
&= \| \bm W^{i} - \bm W^{i+1} \|^{2}_{J} + 2 \| \bm \Sigma^{i+1} - \hat{\bm \Sigma}^{+} \|^{2}_{F}
+ 2 \langle \bm \Lambda^{i+1} - \bm \Lambda^{i}, \bm \Sigma^{i} - \bm \Sigma^{i+1} \rangle.
\end{align}
Hence, next we only need to show $\langle \bm \Lambda^{i+1} - \bm \Lambda^{i}, \bm \Sigma^{i} - \bm \Sigma^{i+1} \rangle \geq 0$.
Now replacing $i$ instead of $i+1$ in \eqref{equ:proof12} and \eqref{equ:proof13} yields
\begin{align}\label{equ:proof21}
( - \bm \Sigma^{i} + \frac{1}{M}\sum_{k=1}^{M}\hat{\bm \Sigma}_{k} - \bm \Lambda^{i})_{jl} \in
\lambda \partial |\bm \Sigma^{i}_{jl}|, ~ j = 1,\ldots,p, l = 1,\ldots,p,  \mbox{and} ~ j \neq l,
\end{align}
and
\begin{align}\label{equ:proof22}
(\bm \Sigma^{i} - \frac{1}{M}\sum_{k=1}^{M}\hat{\bm \Sigma}_{k})_{jj} + \bm \Lambda^{i}_{jj} = 0, ~ j = 1,\ldots,p.
\end{align}
So \eqref{equ:proof12}, \eqref{equ:proof13}, \eqref{equ:proof21} and \eqref{equ:proof22}, together with
the monotonically non-decreasing property of $\partial | \cdot |$, imply
\begin{align}\label{equ:proof23}
\langle \bm \Sigma^{i} - \bm \Sigma^{i+1}, \bm \Lambda^{i+1} - \bm \Lambda^{i} + \bm \Sigma^{i+1} - \bm \Sigma^{i}\rangle \geq 0.
\end{align}
After a simple algebra of \eqref{equ:proof23}, we have
\begin{align*}
\langle \bm \Sigma^{i} - \bm \Sigma^{i+1}, \bm \Lambda^{i+1} - \bm \Lambda^{i} \rangle \geq \| \bm \Sigma^{i+1} - \bm \Sigma^{i} \|^{2}_{F} \geq 0.
\end{align*}
Hence the last two terms on the right hand side of \eqref{equ:proof20} are both non-negative, which proves Lemma 1.
\end{proof}

\begin{proof}{Proof of Theorem 1}{}
According to Lemma 1, we have  \\
\noindent (a) $\| \bm W^{i} - \bm W^{i+1} \|^{2}_{J} \rightarrow 0$, as $i \rightarrow +\infty$; \\
\noindent (b) $\| \bm W^{+} - \bm W^{i} \|^{2}_{J}$ is non-increasing and thus bounded.\\
\noindent The result (a) indicates that $\bm \Sigma^{i} - \bm \Sigma^{i+1} \rightarrow 0$ and
$\bm \Lambda^{i} - \bm \Lambda^{i+1} \rightarrow 0$. Based on \eqref{admm: lambda}, it is easy to see that
$\bm \Phi^{i} - \bm \Sigma^{i} \rightarrow 0$. On the other hand, (b) indicates that $\bm W^{i}$ lies
in a compact region. Accordingly, there exists a subsequence $\bm W^{i_{j}}$
of $\bm W^{i}$ such that $\bm W^{i_{j}} \rightarrow \bm W^{\ast} = (\bm \Lambda^{\ast}, \bm \Sigma^{\ast})$.
In addition, we also have $\bm \Phi^{i_{j}} \rightarrow \bm \Phi^{\ast} \triangleq \bm \Sigma^{\ast}$.
Therefore, $\lim\limits_{i \rightarrow \infty} (\bm \Sigma^{i}, \bm \Phi^{i}, \bm \Lambda^{i}) = (\bm \Sigma^{\ast}, \bm \Phi^{\ast}, \bm \Lambda^{\ast})$.

Next we show that $(\bm \Sigma^{\ast}, \bm \Phi^{\ast}, \bm \Lambda^{\ast})$ is an
optimal solution of \eqref{sigma:equ1}.
By letting $i \rightarrow +\infty$ in \eqref{equ:proof12}, \eqref{equ:proof13} and \eqref{equ:proof6}, we have
\begin{align}\label{equ:proof24}
( - \bm \Sigma^{\ast} + \frac{1}{M}\sum_{k=1}^{M}\hat{\bm \Sigma}_{k} - \bm \Lambda^{\ast})_{jl} \in
\lambda \partial |\bm \Sigma^{\ast}_{jl}|, ~ j = 1,\ldots,p, l = 1,\ldots,p,  \mbox{and} ~ j \neq l,
\end{align}
\begin{align}\label{equ:proof25}
(\bm \Sigma^{\ast} - \frac{1}{M}\sum_{k=1}^{M}\hat{\bm \Sigma}_{k})_{jj} + \bm \Lambda^{\ast}_{jj} = 0, ~ j = 1,\ldots,p,
\end{align}
and
\begin{align}\label{equ:proof26}
\langle \bm \Lambda^{\ast}, \bm \Phi - \bm \Phi^{\ast} \rangle \leq 0,
~ \forall \bm \Phi \succeq \nu \bm I.
\end{align}
\eqref{equ:proof24}, \eqref{equ:proof25} and \eqref{equ:proof26}, together with $\bm \Phi^{\ast} =\bm \Sigma^{\ast}$,
imply that $(\bm \Sigma^{\ast}, \bm \Phi^{\ast}, \bm \Lambda^{\ast})$ is an optimal solution of $\arg \min L(\bm \Sigma, \bm \Phi; \bm \Lambda)$ in \eqref{sigma:optimization}.
Hence, we prove that the sequence produced by Algorithm 1
from any starting point converges to an optimal minimizer of \eqref{sigma:optimization}.
\end{proof}

\begin{proof}{Proof of Lemma 2}{}
The proof is very similar to that of Lemma A.2 in Jiang (2012), so we omit here.
\end{proof}

\begin{proof}{Proof of Lemma 3}{}
We prove this Lemma with the idea of Jiang (2012) by constructing a function $G(\cdot, \cdot)$ via the likelihood function, then decomposing $G(\cdot, \cdot)$ into several parts and considering to bound each part separately.

To simplify the notations, we prove it under the original order without the symbol $\pi_{k}$.
Note that the estimates $\hat{\bm L}$ and $\hat{\bm D}$ based on a sequence of regressions are derived from $\bm \epsilon = \bm L^{-1} \bm X \sim \mathcal{N}(\bm 0, \bm D)$.
The loss functions for the sequence of regressions can be written as the negative log likelihood,
$\sum_{i=1}^{n}[\log \left| \bm D \right| + \tr (\bm x_{i}' \bm L'^{-1} \bm D^{-1} \bm L^{-1} \bm x_{i})]$, up to some constant.
Consequently, adding the penalty terms to the negative log likelihood leads to the following objective function
\begin{align*}
\sum_{i=1}^{n}[\log \left| \bm D \right| + \tr (\bm x_{i}' \bm L'^{-1} \bm D^{-1} \bm L^{-1} \bm x_{i})] + \sum_{j=1}^{p} \eta_{j} \sum_{k<j}|l_{jk}|.
\end{align*}
Denote
\begin{align*}
Q(\bm D, \bm L) = (\log \left| \bm D \right| + \tr (\bm L'^{-1} \bm D^{-1} \bm L^{-1} \bm S) + \sum_{j=1}^{p} \eta_{j} \sum_{k<j}|l_{jk}|.
\end{align*}
Define $G(\Delta_{L}, \Delta_{D}) = Q(\bm D_{0} + \Delta_{D}, \bm L_{0} + \Delta_{L}) - Q(\bm D_{0}, \bm L_{0})$.
Let $\mathcal{A}_{U_{1}} = \{ \Delta_{L}: \| \Delta_{L} \|_{F}^{2} \leq U_{1}^{2} s_{1} \log(p) / n \}$ and
$\mathcal{B}_{U_{2}} = \{ \Delta_{D}: \| \Delta_{D} \|_{F}^{2} \leq U_{2}^{2} p \log(p) / n \}$, where $U_{1}$ and $U_{2}$ are constants.
We will show that for each $\Delta_{L} \in \partial \mathcal{A}_{U_{1}}$ and $\Delta_{D} \in \partial \mathcal{B}_{U_{2}}$, probability $P(G(\Delta_{L}, \Delta_{D}) > 0)$ is tending to 1 as $n \rightarrow \infty$ for sufficiently large $U_{1}$ and $U_{2}$, where $\partial \mathcal{A}_{U_{1}}$ and $\partial \mathcal{B}_{U_{2}}$ are the boundaries of $\mathcal{A}_{U_{1}}$ and $\mathcal{B}_{U_{2}}$, respectively. Additionally, since $G(\Delta_{L}, \Delta_{D}) = 0$ when $\Delta_{L} = 0$ and $\Delta_{D} = 0$, the minimum point of $G(\Delta_{L}, \Delta_{D})$ is achieved when
$\Delta_{L} \in \mathcal{A}_{U_{1}}$ and $\Delta_{D} \in \mathcal{B}_{U_{2}}$.
That is $\| \Delta_{L} \|_{F}^{2} = O_{p} (s_{1} \log(p) / n)$ and $\| \Delta_{D} \|_{F}^{2} = O_{p} (p \log(p) / n)$.

Assume $\| \Delta_{L} \|_{F}^{2} = U_{1}^{2} s_{1} \log(p) / n$ and $\| \Delta_{D} \|_{F}^{2} = U_{2}^{2} p \log(p) / n$. From assumption \eqref{assumption1} and by Lemma 2, without loss of generality,
there exists a constant $h$ such that $0 < 1/h < sv_{p}(\bm L_{0}) \leq sv_{1}(\bm L_{0}) < h < \infty$ and $0 < 1/h < sv_{p}(\bm D_{0}) \leq sv_{1}(\bm D_{0}) < h < \infty$.
Write $\bm D = \bm D_{0} + \Delta_{D}$ and $\bm L = \bm L_{0} + \Delta_{L}$,
then we decompose $G(\Delta_{L}, \Delta_{D})$ into three parts and then consider them separately.
\begin{align*}
G(\Delta_{L}, \Delta_{D}) &= Q(\bm D, \bm L) - Q(\bm D_{0}, \bm L_{0}) \\
&= \log \left| \bm D \right| - \log \left| \bm D_{0} \right| + \tr (\bm L'^{-1} \bm D^{-1} \bm L^{-1} \bm S) - \tr (\bm L'^{-1}_{0} \bm D^{-1}_{0} \bm L^{-1}_{0} \bm S) \\
&~~~+ \sum_{j=1}^{p} \eta_{j} \sum_{k<j}|l_{jk}| - \sum_{j=1}^{p} \eta_{j} \sum_{k<j}|l_{0jk}| \\
&= \log \left| \bm D \right| - \log \left| \bm D_{0} \right| + \tr[(\bm D^{-1} - \bm D_{0}^{-1})\bm D_{0}] - \tr[(\bm D^{-1} - \bm D_{0}^{-1})\bm D_{0}] \\
&~~~+ \tr (\bm L'^{-1} \bm D^{-1} \bm L^{-1} \bm S) - \tr (\bm L'^{-1}_{0} \bm D^{-1}_{0} \bm L^{-1}_{0} \bm S)
+ \sum_{j=1}^{p} \eta_{j} \sum_{k<j}|l_{jk}| - \sum_{j=1}^{p} \eta_{j} \sum_{k<j}|l_{0jk}| \\
&= M_{1} + M_{2} + M_{3},
\end{align*}
where
\begin{align*}
M_{1} &= \log \left| \bm D \right| - \log \left| \bm D_{0} \right| + \tr[(\bm D^{-1} - \bm D_{0}^{-1})\bm D_{0}],\\
M_{2} &= \tr (\bm L'^{-1} \bm D^{-1} \bm L^{-1} \bm S) - \tr (\bm L'^{-1}_{0} \bm D^{-1}_{0} \bm L^{-1}_{0} \bm S) - \tr[(\bm D^{-1} - \bm D_{0}^{-1})\bm D_{0}], \\
M_{3} &= \sum_{j=1}^{p} \eta_{j} \sum_{k<j}|l_{jk}| - \sum_{j=1}^{p} \eta_{j} \sum_{k<j}|l_{0jk}|.
\end{align*}
Based on the proof of Theorem 3.1 in Jiang (2012), we can have $M_{1} \geq 1/8h^4 \|\Delta_{D}\|_{F}^{2}$.
For the second term,
\begin{align*}
M_{2} &= \tr (\bm L'^{-1} \bm D^{-1} \bm L^{-1} \bm S) - \tr (\bm L'^{-1} \bm D^{-1}_{0} \bm L^{-1} \bm S) + \tr (\bm L'^{-1} \bm D^{-1}_{0} \bm L^{-1} \bm S) \\
&~~~ - \tr (\bm L'^{-1}_{0} \bm D^{-1}_{0} \bm L^{-1}_{0} \bm S) - \tr[(\bm D^{-1} - \bm D_{0}^{-1})\bm D_{0}] \\
&= \tr(\bm D^{-1} - \bm D_{0}^{-1})[\bm L^{-1} (\bm S - \bm \Sigma_{0}) \bm L'^{-1}] + \tr \bm D_{0}^{-1} (\bm L^{-1} \bm S \bm L'^{-1} - \bm L^{-1}_{0} \bm S \bm L'^{-1}_{0})    \\
&~~~ + \tr(\bm D^{-1} - \bm D_{0}^{-1})(\bm L^{-1} \bm \Sigma_{0} \bm L'^{-1} - \bm D_{0})  \\
&= \tr(\bm D^{-1} - \bm D_{0}^{-1})[\bm L^{-1} (\bm S - \bm \Sigma_{0}) \bm L'^{-1}] + \tr [\bm D_{0}^{-1} (\bm L^{-1} (\bm S - \bm \Sigma_{0}) \bm L'^{-1} - \bm L^{-1}_{0} (\bm S - \bm \Sigma_{0}) \bm L'^{-1}_{0})]  \\
&~~~ + \tr [\bm D_{0}^{-1} (\bm L^{-1} \bm \Sigma_{0} \bm L'^{-1} - \bm L^{-1}_{0} \bm \Sigma_{0} \bm L'^{-1}_{0})] + \tr(\bm D^{-1} - \bm D_{0}^{-1})(\bm L^{-1} \bm \Sigma_{0} \bm L'^{-1} - \bm D_{0})  \\
&= \tr(\bm D^{-1} - \bm D_{0}^{-1})[\bm L^{-1} (\bm S - \bm \Sigma_{0}) \bm L'^{-1}] + \tr [\bm D_{0}^{-1} (\bm L^{-1} (\bm S - \bm \Sigma_{0}) \bm L'^{-1} - \bm L^{-1}_{0} (\bm S - \bm \Sigma_{0}) \bm L'^{-1}_{0})]  \\
&~~~ + \tr [\bm D^{-1} (\bm L^{-1} \bm \Sigma_{0} \bm L'^{-1} - \bm L^{-1}_{0} \bm \Sigma_{0} \bm L'^{-1}_{0})] \\
&= M_{2}^{(1)} + M_{2}^{(2)} + M_{2}^{(3)},
\end{align*}
where the fourth equality uses $\bm L^{-1}_{0} \bm \Sigma_{0} \bm L'^{-1}_{0} = \bm L^{-1}_{0} (\bm L_{0} \bm D_{0} \bm L_{0}^{'}) \bm L'^{-1}_{0} = \bm D_{0}$. The notations $M_{2}^{(1)}$, $M_{2}^{(2)}$ and $M_{2}^{(3)}$ are defined in the following
\begin{align*}
M_{2}^{(1)} &= \tr(\bm D^{-1} - \bm D_{0}^{-1})[\bm L^{-1} (\bm S - \bm \Sigma_{0}) \bm L'^{-1}],\\
M_{2}^{(2)} &= \tr [\bm D_{0}^{-1} (\bm L^{-1} (\bm S - \bm \Sigma_{0}) \bm L'^{-1} - \bm L^{-1}_{0} (\bm S - \bm \Sigma_{0}) \bm L'^{-1}_{0})], \\
M_{2}^{(3)} &= \tr [\bm D^{-1} (\bm L^{-1} \bm \Sigma_{0} \bm L'^{-1} - \bm L^{-1}_{0} \bm \Sigma_{0} \bm L'^{-1}_{0})].
\end{align*}
Based on the proof of Theorem 3.1 in Jiang (2012), for any $\epsilon > 0$, there exists $V_{1} > 0$ and $V_{2} > 0$ such that
\begin{align*}
|M_{2}^{(1)}| \leq V_{1} \sqrt{p \log (p) / n} \|\Delta_{D}\|_{F}
\end{align*}
and
\begin{align*}
M_{2}^{(3)} - |M_{2}^{(2)}| > 1 /2 h^4 \| \Delta_{L} \|_{F}^2 - V_{2} \sqrt{\log (p) / n} \sum_{(j, k) \in Z^c} | l_{jk} | - V_{2} h \sqrt{s_{1} \log (p) / n} \| \Delta_{L} \|_{F},
\end{align*}
where $Z = \{(j, k): k < j, l_{0jk} \neq 0 \}$, and $l_{0jk}$ represents the element $(j, k)$ of the matrix $\bm L_{0}$.
Next, for the penalty term,
\begin{align*}
M_{3} = \sum_{j=1}^{p} \eta_{j} \sum_{(j, k) \in Z^c}|l_{jk}| + \sum_{j=1}^{p} \eta_{j}
\sum_{(j, k) \in Z} (|l_{jk}| - |l_{0jk}|) = M_{3}^{(1)} + M_{3}^{(2)},
\end{align*}
where
\begin{align*}
M_{3}^{(1)} = \sum_{j=1}^{p} \eta_{j} \sum_{(j, k) \in Z^c}|l_{jk}|,
\end{align*}
and
\begin{align*}
| M_{3}^{(2)} | = | \sum_{j=1}^{p} \eta_{j} \sum_{(j, k) \in Z}(|l_{jk}| - |l_{0jk}|) |
&\leq \sum_{j=1}^{p} \eta_{j} \sum_{(j, k) \in Z}| (|l_{jk}| - |l_{0jk}|) |        \\
&\leq \sum_{j=1}^{p} \eta_{j} \sum_{(j, k) \in Z}| l_{jk} - l_{0jk} |   \\
&\leq \sum_{j=1}^{p} \eta_{j} \sqrt{s_{1}} \| \Delta_{L} \|_{F},
\end{align*}
where the last inequality uses the fact that $(a_{1} + a_{2} + \cdots + a_{m})^2 \leq m (a_{1}^2 + a_{2}^2 + \cdots + a_{m}^2)$.
Combine all the terms above together, with probability greater than $1 - 2 \epsilon$, we have
\begin{align*}
&~~~ | G(\Delta_{L}, \Delta_{D}) | \\
&\geq M_{1} - |M_{2}^{(1)}| + M_{2}^{(3)} - |M_{2}^{(2)}| + M_{3}^{(1)} - |M_{3}^{(2)}|   \\
&\geq 1/8h^4 \|\Delta_{D}\|_{F}^{2} - V_{1} \sqrt{p \log (p) / n} \|\Delta_{D}\|_{F} + 1 /2 h^4 \| \Delta_{L} \|_{F}^2 - V_{2} \sqrt{\log (p) / n} \sum_{(j, k) \in Z^c} | l_{jk} |   \\
&~~~ - V_{2} h \sqrt{s_{1} \log (p) / n} \| \Delta_{L} \|_{F} + \sum_{j=1}^{p} \eta_{j} \sum_{(j, k) \in Z^c}|l_{jk}| - \sum_{j=1}^{p} \eta_{j} \sqrt{s_{1}} \| \Delta_{L} \|_{F}    \\
&= \frac{U_{2}^{2}}{8h^4} p \log(p) / n - V_{1} U_{2} p \log(p) / n + \frac{U_{1}^{2}}{2h^{4}} s_{1} \log(p) / n - V_{2} \sqrt{\log (p) / n} \sum_{(j, k) \in Z^c} | l_{jk} |  \\
&~~~ - V_{2} U_{1} h s_{1} \log(p) / n + \sum_{j=1}^{p} \eta_{j} \sum_{(j, k) \in Z^c}|l_{jk}| - s_{1} U_{1} \sqrt{\log (p) / n} \sum_{j=1}^{p} \eta_{j} \\
&= \frac{U_{2} p \log(p)}{n} (\frac{U_{2}}{8h^4} - V_{1}) + \frac{U_{1} s_{1} p \log(p)}{n} (\frac{U_{1}}{2h^{4}} - \frac{\sum_{j=1}^{p} \eta_{j}}{\sqrt{\log (p) / n}} - V_{2} h )     \\
&~~~ + \sum_{(j, k) \in Z^c}|l_{jk}| (\sum_{j=1}^{p} \eta_{j} - V_{2} \sqrt{\log (p) / n}).
\end{align*}
Here $V_{1}$ and $V_{2}$ are only related to the sample size $n$ and $\epsilon$. Assume $\sum_{j=1}^{p} \eta_{j} = K(\log (p) / n)$ where $K > V_{2}$ and choose $U_{1} > 2h^4(K + hV_{2})$, $U_{2} > 8h^4 V_{1}$, then $G(\Delta_{L}, \Delta_{D}) > 0$. This establishes the lemma.
\end{proof}

\begin{proof}{Proof of Theorem 2}{}
Based on the proof of Theorem 3.2 in Jiang (2012),
we can have
\begin{align*}
\| \hat{\bm \Sigma}_{\pi_{k}} - \bm \Sigma_{0\pi_{k}} \|_{F}^2
&= O_{p}(\| \hat{\bm L}_{\pi_{k}} - \bm L_{0\pi_{k}} \|^2_{F}) + O_{p}(\| \hat{\bm D}_{\pi_{k}} - \bm D_{0\pi_{k}} \|^2_{F})   \\
&= O_{p}(s_{1}\log (p) / n) + O_{p}(p \log (p) / n) \\
&= O_{p}((s_{1} + p) \log (p) / n),
\end{align*}
where the second inequality is provided by the proof of Lemma 3. Then
\begin{align*}
\| \hat{\bm \Sigma}_{k} - \bm \Sigma_{0} \|_{F}^2 &= \|  \bm P_{\pi_{k}} \hat{\bm \Sigma}_{\pi_{k}} \bm P_{\pi_{k}}' - \bm P_{\pi_{k}} \bm \Sigma_{0\pi_{k}} \bm P_{\pi_{k}}' \|_{F}^2  \\
&= \| \bm P_{\pi_{k}} (\hat{\bm \Sigma}_{\pi_{k}} - \bm \Sigma_{0\pi_{k}}) \bm P_{\pi_{k}}' \|_{F}^2   \\
&= \| \hat{\bm \Sigma}_{\pi_{k}} - \bm \Sigma_{0\pi_{k}} \|_{F}^2      \\
&= O_{p}((s_{1} + p) \log (p) / n),
\end{align*}
where the third equality uses the fact that the Frobenius norm of a matrix is invariant on the permutation matrix.

As $\bm \Sigma_{0}$ is positive definite, there exists $\epsilon > 0$ such that $\epsilon < \lambda_{min}(\bm \Sigma_{0})$, where $\lambda_{min}(\bm \Sigma_{0})$ is the smallest eigenvalue of $\bm \Sigma_{0}$.
By introducing $\bm \Delta = \bm \Sigma - \bm \Sigma_{0}$, the expression of \eqref{sigma:equ1} can be rewritten in terms of $\bm \Delta$ as
\begin{align*}
\hat{\bm \Delta} = \arg \min_{\bm \Delta=\bm \Delta', \bm \Delta + \bm \Sigma_{0} \succeq \epsilon \bm I} \frac{1}{2M} \sum_{k=1}^{M} \| \bm \Delta + \bm \Sigma_{0} -  \hat{\bm \Sigma}_{k} \|_{F}^{2} + \lambda | \bm \Delta + \bm \Sigma_{0} |_{1}~~(\triangleq \mathcal{F}(\bm \Delta)).
\end{align*}
Note that it is easy to see $\hat{\bm \Delta} = \hat{\bm \Sigma}^{+} - \bm \Sigma_{0}$.
Now consider $\bm \Delta \in \{ \bm \Delta: \bm \Delta=\bm \Delta', \bm \Delta + \bm \Sigma_{0} \succeq \epsilon \bm I, \|\bm \Delta\|_{F} = 5 \lambda \sqrt{s_{0} + p} \}$.
Define the active set of $\bm \Sigma_{0}$ as $A_{0} = \{(i,j): \sigma_{ij}^0 \neq 0, i \neq j\}$,
and $\bm B_{A_{0}} = (b_{ij} \cdot \bm I_{\{(i,j) \in A_{0}\}})_{1 \leq i, j \leq p}$.
Let $A_{0}^{c}$ be the complement set of $A_{0}$.
Denote the element $(i,j)$ of matrix $\bm \Delta$ by $\bm \Delta_{ij}$.
Under the probability event $\{ |\hat{\sigma}_{ij}^k - \sigma_{ij}^{0}| \leq \lambda \}$ where $\hat{\bm \Sigma}_{k} = (\hat{\sigma}_{ij}^k)_{p \times p}$, we have

\begin{align*}
\mathcal{F}(\bm \Delta) - \mathcal{F}(0) &= \frac{1}{2M} \sum_{k=1}^{M} \| \bm \Delta + \bm \Sigma_{0} -  \hat{\bm \Sigma}_{k} \|_{F}^{2} - \frac{1}{2M} \sum_{k=1}^{M} \| \bm \Sigma_{0} -  \hat{\bm \Sigma}_{k} \|_{F}^{2} + \lambda | \bm \Delta + \bm \Sigma_{0} |_{1} - \lambda |\bm \Sigma_{0} |_{1} \\
&= \frac{1}{2}\| \bm \Delta \|_{F}^{2} + \frac{1}{M} \sum_{k=1}^{M} <\bm \Delta, \bm \Sigma_{0} -  \hat{\bm \Sigma}_{k}> + \lambda |\bm \Delta_{A_{0}^c}|_{1} + \lambda (|\bm \Delta_{A_{0}} + (\bm \Sigma_{0})_{A_{0}}|_{1} \\
&~~~~~~ - |(\bm \Sigma_{0})_{A_{0}}|_{1}) \\
&\geq \frac{1}{2}\| \bm \Delta \|_{F}^{2} - \lambda (|\bm \Delta|_{1} + \sum_{i} \bm \Delta_{ii}) + \lambda |\bm \Delta_{A_{0}^c}|_{1} - \lambda |\bm \Delta_{A_{0}}|_{1}    \\
&\geq \frac{1}{2}\| \bm \Delta \|_{F}^{2} - 2 \lambda (|\bm \Delta_{A_{0}}|_{1} + \sum_{i} \bm \Delta_{ii})  \\
&\geq \frac{1}{2}\| \bm \Delta \|_{F}^{2} - 2 \lambda \sqrt{s_{0} + p} \| \bm \Delta \|_{F}   \\
&= \frac{5}{2} \lambda^2 (s_{0} + p) \\
&> 0.
\end{align*}
Note that $\hat{\bm \Delta}$ is also the optimal solution to the convex optimization problem
\begin{align*}
\hat{\bm \Delta} = \arg \min_{\bm \Delta=\bm \Delta', \bm \Delta + \bm \Sigma_{0} \succeq \epsilon \bm I} \mathcal{F}(\bm \Delta) - \mathcal{F}(0).
\end{align*}
The rest of proof is the same as that of Theorem 2 in Xue, Ma, \& Zou (2012), so we omit them.
\end{proof}
\end{appendix}

\begin{table}
\footnotesize
\begin{center}
\caption{The averages and standard errors of estimates for $\textbf{Model}$ 3.}
\label{sigma:table_model3}
\begin{tabular}{rrrrrrrrrrrrrrr}
\hline\hline
&  &F &EN &$L_{1}$ &MAE &FSL (\%) \\\hline
\multirow {9}*{$p = 30$}
&$\bm S$&4.39 (0.04)      & 12.46 (0.09)    & 5.08 (0.06)     & 3.48 (0.02)     & 83.96 (0.01) \\
&BIC    &3.28 (0.03)      & 7.28 (0.10)     & 2.93 (0.04)     & 1.76 (0.01)     & 52.71 (0.51) \\
&BPA    &3.34 (0.03)      & \textbf{5.87 (0.09)} & 2.60 (0.05)    & 1.52 (0.02)     & 46.59 (0.78) \\
&BT     &4.73 (0.01)      & 7.77 (0.04)     & 2.14 (0.01)     & 1.85 (0.00)     & 6.94 (0.10) \\
&BL     &3.28 (0.05)      & -               & 2.24 (0.05)     & \textbf{1.14 (0.01)} & \textbf{6.62 (0.15)} \\
&XMZ    &3.33 (0.04)      & 10.53 (0.14)    & \textbf{1.87 (0.01)} & 1.23 (0.01) & 7.00 (0.15) \\
&IB     &\textbf{2.99 (0.04)}      & -               & 2.24 (0.04)    & 1.26 (0.03)    & 9.37 (0.44) \\
&RLZ    &4.37 (0.01)      & 16.75 (0.11)    & 2.28 (0.02) & 1.68 (0.00) & 15.55 (0.00) \\
&Proposed   &\textbf{3.22 (0.04)}  & \textbf{6.96 (0.10)} & \textbf{1.89 (0.02)} & \textbf{1.21 (0.01)} & \textbf{6.84 (0.13)} \\
\midrule
\hline
\multirow {9}*{$p = 50$}
&$\bm S$&7.25 (0.05)      & -               & 8.28 (0.08)     & 5.75 (0.03)     & 90.19 (0.01) \\
&BIC    &4.91 (0.03)      & 16.00 (0.21)    & 3.85 (0.08)     & 1.98 (0.02)     & 43.07 (0.50) \\
&BPA    &4.89 (0.03)      & \textbf{12.92 (0.16)} & 3.55 (0.08)  & 1.83 (0.02)     & 40.90 (0.65) \\
&BT     &6.18 (0.06)      & 15.58 (0.26)    & 2.45 (0.02)     & 1.96 (0.01)     & 11.21 (0.47) \\
&BL     &4.65 (0.05)      & -               & 2.46 (0.05)     & \textbf{1.27 (0.01)} & \textbf{4.79 (0.06)} \\
&XMZ    &4.63 (0.06)      & 20.40 (0.32)    & \textbf{1.97 (0.01)} & 1.36 (0.01) & 5.18 (0.07) \\
&IB     &\textbf{4.35 (0.04)}      & -               & 2.53 (0.05)    & \textbf{1.25 (0.02)}    & 6.48 (0.32) \\
&RLZ    &6.03 (0.01)      & 31.67 (0.17)    & 2.43 (0.01) & 1.90 (0.00)    & 11.36 (0.00) \\
&Proposed   &\textbf{4.60 (0.03)}  & \textbf{13.74 (0.11)} & \textbf{2.00 (0.01)} & 1.36 (0.01) & \textbf{4.39 (0.07)} \\
\midrule
\hline
\multirow {9}*{$p = 100$}
&$\bm S$&14.41 (0.05)     & -               & 16.18 (0.10)    & 11.43 (0.03)    & 95.01 (0.00) \\
&BIC    &7.59 (0.03)      & 42.88 (0.52)    & 5.36 (0.10)     & 2.26 (0.01)     & 32.65 (0.37) \\
&BPA    &7.14 (0.03)      & 35.37 (0.43)    & 5.12 (0.12)   & 2.19 (0.02)     & 33.30 (0.41) \\
&BT     &8.62 (0.14)      & \textbf{28.78 (0.37)} & \textbf{2.36 (0.02)} & 1.89 (0.02) & 3.91 (0.30) \\
&BL     &7.18 (0.05)      & -               & 2.58 (0.05) & \textbf{1.41 (0.01)} & 2.83 (0.02) \\
&XMZ    &14.40 (0.05)     & 364.17 (0.16)   & 16.17 (0.10)    & 11.42 (0.03)    & 94.96 (0.00) \\
&IB     &\textbf{7.11 (0.04)}      & -               & 2.80 (0.05)    & \textbf{1.21 (0.01)}    & \textbf{2.61 (0.07)} \\
&RLZ    &8.57 (0.01)      & 64.68 (0.22)    & 2.50 (0.01) & 1.92 (0.00)    & 5.72 (0.00) \\
&Proposed   &\textbf{7.06 (0.03)}  & \textbf{31.33 (0.15)} & \textbf{2.10 (0.01)} & 1.49 (0.00) & \textbf{2.39 (0.02)} \\
\midrule
\hline
\end{tabular}
\end{center}
\end{table}

\begin{table}
\footnotesize
\begin{center}
\caption{The averages and standard errors of estimates for $\textbf{Model}$ 4.}
\label{sigma:table_model4}
\begin{tabular}{rrrrrrrrrrrrr}
\hline\hline
&  &F &EN &$L_{1}$ &MAE &FSL (\%) \\\hline
\multirow {9}*{$p = 30$}
&$\bm S$&4.41 (0.05)      & 12.48 (0.08)    & 5.16 (0.07)     & 3.49 (0.03)     & 46.66 (0.01) \\
&BIC    &3.35 (0.03)      & 5.22 (0.08)     & 2.89 (0.04)     & 1.90 (0.01)     & 43.03 (0.33) \\
&BPA    &\textbf{3.13 (0.03)}      & \textbf{4.48 (0.08)} & 2.79 (0.04)    & 1.76 (0.01)     & 42.66 (0.35) \\
&BT     &4.69 (0.01)      & 5.36 (0.04)     & 2.50 (0.01)     & 2.19 (0.00)     & 43.79 (0.08) \\
&BL     &3.47 (0.05)      & -               & 2.69 (0.06)     & \textbf{1.58 (0.02)} & 43.86 (0.20) \\
&XMZ    &3.51 (0.03)      & 5.86 (0.09)     & \textbf{2.24 (0.01)} & 1.64 (0.01) & 42.31 (0.19) \\
&IB     &\textbf{3.08 (0.04)}      & 7.91 (0.21)   & 2.62 (0.05)    & \textbf{1.57 (0.02)}    & \textbf{41.28 (0.39)} \\
&RLZ    &4.55 (0.01)      & 11.57 (0.10)    & 2.67 (0.02) & 2.15 (0.00)    & 50.67 (0.00) \\
&Proposed   &3.48 (0.03)  & \textbf{4.06 (0.05)} & \textbf{2.27 (0.01)} & 1.66 (0.01) & \textbf{42.06 (0.16)} \\
\midrule
\hline
\multirow {9}*{$p = 50$}
&$\bm S$&7.24 (0.05)      & -               & 8.24 (0.08)     & 5.74 (0.03)     & 65.58 (0.01) \\
&BIC    &4.59 (0.03)      & 11.16 (0.21)    & 3.95 (0.08)     & 2.17 (0.01)     & 42.10 (0.24) \\
&BPA    &\textbf{4.40 (0.03)}      & \textbf{9.22 (0.15)} & 3.72 (0.08)   & 2.07 (0.01)     & 41.28 (0.30) \\
&BT     &6.09 (0.04)      & 10.41 (0.21)    & 2.71 (0.01)     & 2.27 (0.01)     & 30.31 (0.14) \\
&BL     &4.71 (0.04)      & -               & 2.79 (0.06)     & \textbf{1.68 (0.01)} & 29.50 (0.07) \\
&XMZ    &4.72 (0.04)      & 10.96 (0.19)    & \textbf{2.37 (0.01)} & 1.75 (0.01) & 28.78 (0.08) \\
&IB     &\textbf{4.21 (0.04)}      & 14.23 (0.31)   & 2.86 (0.04)    & \textbf{1.62 (0.02)}    & \textbf{27.98 (0.20)} \\
&RLZ    &5.94 (0.01)      & 19.21 (0.11)    & 2.77 (0.01) & 2.24 (0.00)    & 33.60 (0.00) \\
&Proposed   &4.74 (0.03)  & \textbf{7.29 (0.06)} & \textbf{2.41 (0.01)} & 1.78 (0.01) & \textbf{28.56 (0.07)} \\
\midrule
\hline
\multirow {9}*{$p = 100$}
&$\bm S$&14.41 (0.08)     & -               & 16.10 (0.13)    & 11.44 (0.04)    & 81.85 (0.00) \\
&BIC    &6.92 (0.02)      & 30.15 (0.51)    & 5.39 (0.10) & 2.46 (0.01)     & 33.63 (0.25) \\
&BPA    &\textbf{6.80 (0.03)}      & 23.53 (0.35)    & 5.39 (0.13)   & 2.43 (0.01)     & 33.92 (0.29) \\
&BT     &8.51 (0.11)      & \textbf{20.62 (0.33)} & \textbf{2.73 (0.02)} & 2.24 (0.02) & 16.43 (0.17) \\
&BL     &7.20 (0.04)      & -               & 3.03 (0.06) & \textbf{1.82 (0.01)} &16.08 (0.02) \\
&XMZ    &14.40 (0.08)     & 369.42 (0.22)   & 16.09 (0.13)    & 11.43 (0.04) & 81.82 (0.01) \\
&IB    &\textbf{6.76 (0.05)}        & 27.14 (0.39)   & 3.23 (0.04)    & \textbf{1.63 (0.01)}    & \textbf{15.12 (0.06)} \\
&RLZ    &8.54 (0.01)      & 40.19 (0.21)    & 2.89 (0.01) & 2.31 (0.00) & 18.44 (0.00) \\
&Proposed   &7.12 (0.03)  & \textbf{16.39 (0.08)} & \textbf{2.49 (0.01)} & 1.90 (0.00) & \textbf{15.64 (0.02)} \\
\midrule
\hline
\end{tabular}
\end{center}
\end{table}

\begin{table}
\footnotesize
\begin{center}
\caption{The averages and standard errors of estimates for $\textbf{Model}$ 5.}
\label{sigma:table_model5}
\begin{tabular}{rrrrrrrrrrrrr}
\hline\hline
&  &F &EN &$L_{1}$ &MAE &FSL (\%) \\\hline
\multirow {9}*{$p = 30$}
&$\bm S$& 4.38 (0.03) & 12.52 (0.08) & 4.90 (0.05) & 3.47 (0.02) & 83.72 (0.01) \\
&BIC    & \textbf{2.73 (0.02)} & \textbf{6.80 (0.12)} & 2.76 (0.07) & 1.34 (0.02) & 38.29 (0.51) \\
&BPA    & 2.76 (0.03) & \textbf{8.20 (0.13)} & 3.06 (0.09) & 1.46 (0.02) & 55.20 (0.78) \\
&BT     & 3.49 (0.01) & 10.15 (0.15) & 2.06 (0.01) & 1.21 (0.01) & \textbf{10.23 (0.09)} \\
&BL     & 2.88 (0.02) & -           & 1.92 (0.03) & 0.96 (0.01) & 11.65 (0.08) \\
&XMZ    & 2.78 (0.02) & 19.72 (0.28) & \textbf{1.78 (0.02)} & \textbf{0.94 (0.01)} & 11.45 (0.09) \\
&IB     & 2.90 (0.03) & -           & 2.01 (0.04) & 1.08 (0.02) & 14.97 (0.50) \\
&RLZ    & 3.21 (0.01) & 23.74 (0.19) & 2.18 (0.01) & 1.22 (0.01) & 18.44 (0.01) \\
&Proposed & \textbf{2.48 (0.02)} & 11.80 (0.29) & \textbf{1.62 (0.02)} & \textbf{0.86 (0.01)} & \textbf{10.73 (0.11)} \\
\midrule
\hline
\multirow {9}*{$p = 50$}
&$\bm S$& 7.25 (0.03) & -            & 8.03 (0.06) & 5.75 (0.02) & 84.43 (0.01) \\
&BIC    & 3.92 (0.02) & \textbf{21.74 (0.41)} & 4.09 (0.11) & 1.71 (0.01) & 30.76 (0.37) \\
&BPA    & \textbf{3.75 (0.02)} & 22.41 (0.28) & 3.78 (0.09) & 1.76 (0.02) & 47.28 (0.72) \\
&BT     & 4.46 (0.01) & \textbf{20.83 (0.21)} & 2.58 (0.01) & 1.48 (0.01) & \textbf{13.25 (0.13)} \\
&BL     & 3.79 (0.01) & -            & 2.44 (0.02) & 1.26 (0.01) & 13.37 (0.02) \\
&XMZ    & 3.68 (0.01) & 38.03 (0.23) & \textbf{2.38 (0.01)} & \textbf{1.24 (0.01)} & 13.49 (0.01) \\
&IB     & 4.07 (0.01) & -            & 2.63 (0.02) & 1.44 (0.01) & 14.57 (0.05) \\
&RLZ    & 3.90 (0.01) & 37.04 (0.25) & 2.48 (0.02) & 1.44 (0.01) & 16.96 (0.01) \\
&Proposed & \textbf{3.63 (0.01)} & 32.08 (0.26) & \textbf{2.31 (0.01)} & \textbf{1.23 (0.01)} & \textbf{13.21 (0.04)} \\
\midrule
\hline
\multirow {9}*{$p = 100$}
&$\bm S$& 14.26 (0.04) & -           & 15.56 (0.09) & 11.32 (0.03) & 84.32 (0.01) \\
&BIC    & 5.87 (0.02) & 43.75 (0.47) & 5.38 (0.11) & 2.27 (0.01) & 24.49 (0.22) \\
&BPA    & 5.77 (0.02) & 37.55 (0.28) & 5.67 (0.13) & 2.39 (0.02) & 38.77 (0.45) \\
&BT     & 6.76 (0.02) & \textbf{31.42 (0.30)} & 3.37 (0.02) & 2.12 (0.01) & 15.99 (0.26) \\
&BL     & \textbf{5.38 (0.02)} & -            & \textbf{3.09 (0.02)} & \textbf{1.87 (0.01)} & \textbf{14.69 (0.01)} \\
&XMZ    & 14.26 (0.04) & 368.48 (0.19) & 15.55 (0.09) & 11.31 (0.03) & 84.28 (0.01) \\
&IB     & 6.07 (0.01) & -            & 3.35 (0.02) & 2.05 (0.01) & 15.62 (0.02) \\
&RLZ    & 5.51 (0.01) & 38.82 (0.12) & 3.23 (0.02) & 1.94 (0.01) & 16.22 (0.01) \\
&Proposed & \textbf{5.23 (0.01)} & \textbf{33.60 (0.09)} & \textbf{3.06 (0.01)} & \textbf{1.76 (0.01)} & \textbf{14.62 (0.01)} \\
\midrule
\hline
\end{tabular}
\end{center}
\end{table}

\end{document}